\documentclass[12pt,english]{article}
\usepackage[T1]{fontenc}
\usepackage[latin9]{inputenc}
\usepackage[a4paper]{geometry}
\usepackage{graphicx}
\usepackage{color}
\usepackage{indentfirst}
\usepackage{enumerate}
\usepackage{amsmath}
\usepackage{amssymb}
\usepackage{tikz}
\usepackage{leftindex}
\usepackage{comment}
\usepackage{subcaption}
\usepackage[arrow]{xy}
\usepackage[linktocpage=true]{hyperref}
\usepackage{tocloft} 
\usepackage{babel}
\usepackage[T1]{fontenc}

\newtheorem{thm}{Theorem}[section]
\newtheorem{definition}[thm]{Definition}
\newtheorem{remark}[thm]{Remark}
\newtheorem{proposition}[thm]{Proposition}
\newtheorem{corollary}[thm]{Corollary}
\newtheorem{example}[thm]{Example}
\newtheorem{lemma}[thm]{Lemma}

\usepackage{authblk}
\makeatletter
\def\thanks#1{\protected@xdef\@thanks{\@thanks
		\protect\footnotetext{#1}}}
\makeatother

\title{Jacobi Hamiltonian Integrators}
\author{Ad\'erito Ara\'ujo\thanks{alma@mat.uc.pt}}
\author{Gon\c calo Inoc\^encio Oliveira\thanks{g.inoc.oliveira@gmail.com}}
\author{Jo\~ao Nuno Mestre\thanks{jnmestre@proton.me}}
\affil{CMUC, University of Coimbra, Department of Mathematics, Portugal.}

\date{}
\begin{document}
	
		\maketitle
	\begin{abstract}
		We develop a method of constructing structure-preserving integrators for Hamiltonian systems in Jacobi manifolds. Hamiltonian mechanics, rooted in symplectic and Poisson geometry, has long provided a foundation for modeling conservative systems in classical physics. Jacobi manifolds, generalizing both contact and Poisson manifolds, extend this theory and are suitable for incorporating time-dependent, dissipative and thermodynamic phenomena.
		Building on recent advances in geometric integrators - specifically Poisson Hamiltonian Integrators (PHI), which preserve key features of Poisson systems - we propose a construction of Jacobi Hamiltonian Integrators. Our approach explores the correspondence between Jacobi and homogeneous Poisson manifolds, with the aim of extending the PHI techniques while ensuring preservation of the homogeneity structure.
		This work develops the theoretical tools required for this generalization and outlines a numerical integration technique compatible with Jacobi dynamics. { By focusing on the homogeneous Poisson perspective instead of direct contact realizations, we establish a clear pathway for constructing structure-preserving integrators for time-dependent and dissipative systems that are embedded in the Jacobi framework.}
	\end{abstract}
	
	\tableofcontents
	
	
	\section{Introduction}
	
	Symplectic and Poisson geometry have long played a foundational role in physics, particularly through Hamilton's equations, which provide a natural framework for conservative mechanical systems in classical mechanics.
	In order to extend Hamiltonian dynamics to be able describe systems interacting with their environment, it is useful to step out of the symplectic and Poisson frameworks.
	
	Contact geometry offers a direct generalization of Hamiltonian mechanics, which naturally models dissipative and thermodynamic processes \cite{MRUGALA1991109}. Recent developments have underscored the growing relevance of contact Hamiltonian dynamics in thermodynamics \cite{Bravetti_19} and partial differential equations \cite{Kushner_Lychagin_Rubtsov_2006}. Jacobi geometry extends contact geometry in a way similar to how Poisson geometry extends symplectic geometry: by allowing for degeneracy. Jacobi manifolds, introduced by Kirillov \cite{Kirillov} and Lichnerowicz \cite{LICHNEROWICZ1977}, provide a powerful framework for extending both contact and Poisson manifolds, enabling the treatment of time-dependent and dissipative dynamics.

	Jacobi and contact geometry can be interpreted as extensions, or analogues of, respectively, Poisson and symplectic geometry. But the parallel is in fact much stronger: there is a 1-to-1 identification of Jacobi manifolds (resp. contact) and Poisson manifolds (resp. symplectic) which are homogeneous, meaning that they are equipped with a compatible scaling symmetry (a free and proper action of the non-zero reals). This identification, given by Poissonization of Jacobi geometry \cite{LICHNEROWICZ1977} (resp. symplectization of contact geometry) induces an equivalence of categories \cite{Monica_Blaga_2020}, and it has proven to be a powerful point of view, permitting the study of Jacobi and contact structures with the tools of Poisson and symplectic geometry (e.g. in \cite{DLM1991,BruceGrabowskaGrabowski2017,Monica_Blaga_2020,CMS2022,Grabowska_2022,Jonas2023}). It brings us full-circle back to the symplectic and Poisson framework that we had initially stepped out of, at the (small) price of having to care for an added homogeneity structure, that must be preserved by all results and constructions that we wish to use. That is the strategy of the present paper, for the construction of structure preserving integrators for Hamiltonian systems in Jacobi manifolds.
	
	Recently, a numerical method for solving Poisson Hamiltonian systems was developed \cite{COSSERAT2023,cosserat2023numerical}, called a Poisson Hamiltonian Integrator (PHI). This method lifts
	the system to a local symplectic groupoid, applies Hamilton-Jacobi techniques via Lagrangian bisections, and projects back to the Poisson manifold at each step. This method has proven to be very effective for this type of system, since by preserving the Poisson structure and the Hamiltonian (up to a certain order), it preserves crucial properties of the geometry such as the symplectic foliation and Casimir functions.
	
	In this paper, our aim is to establish all the tools to construct a Jacobi Hamiltonian Integrator using a similar strategy as that for PHI. To do so, we have two alternatives: using the analogous technique for Jacobi manifolds (using the constructions of contact realizations and local contact groupoids of \cite{CMS2022}); or using the correspondence between Jacobi and homogeneous Poisson and proving that each of the tools and results involved the PHI technique can be made to preserve the homogeneity. In this work, we proceed with the second approach and obtain a constructive and geometrically natural method to produce Jacobi Hamiltonian integrators.
	
	In this approach we used versions of the Darboux-Weinstein and of the Weinstein tubular neighborhood Theorems for homogeneous symplectic manifolds and homogeneous Lagrangians submanifolds (Theorems \ref{thm: Darboux-Weinstein} and \ref{thm-WLtub}). We make no claim of originality of these results, which are known and proved in the equivalent context of contact geometry (e.g. in \cite{Loose, Geiges}). We do provide proofs for them that are straightforward generalizations of the usual symplectic versions, by carefully checking compatibility with homogeneity, because they were more easily applied to our constructions. To the best of our knowledge we could not find such proofs in the literature, although there are related ones: for homogeneous versions of the Darboux theorem \cite{DLM1991, Jonas2023, GG2025}, and of the Poincar\'e Lemma \cite{GG2025}, for example; we believe these may be of independent interest.
	
	
	\subsubsection*{Structure of the paper}
	
	In \textbf{Section \ref{anexo: homogeneous}} we recall the background on homogeneous versions of manifolds, smooth maps, differential forms and multivector fields, and submanifolds.
	
	In \textbf{Section \ref{sec: poissonization}} we briefly describe Poisson, Jacobi, and contact manifolds, and the Poissonization procedure that lets us interpret Jacobi geometry as homogeneous Poisson geometry.
	
	\textbf{Section \ref{sec-homogeneous-tools}} introduces some tools from homogeneous Poisson and symplectic geometry: in the first part we describe homogeneous Lagrangian submanifolds, and homogeneous versions of the Darboux-Weinstein and the Weinstein Lagrangian neighborhood Theorems; in the second part we describe a construction of homogeneous symplectic bi-realizations, in terms of homogeneous Poisson sprays.
	
	\textbf{Section \ref{sec-smoothfamilies}} concerns smooth families of homogeneous Lagrangian submanifolds, and specifically of homogeneous Lagrangian bisections of bi-realizations. These will be essential in our construction of Jacobi Hamiltonian integrators. { A numerical example of a dissipative contact system is included, illustrating the ability of the Jacobi Hamiltonian integrator to accurately capture dissipative dynamics while preserving the underlying contact structure.}
	
	In \textbf{Section \ref{sec: JHI}}, using the tools from the tools from Sections 3, 4 and 5, we proceed to the construction of homogeneous Poisson-Hamilton integrators, we define Jacobi-Hamilton integrators, and we establish a relation between both.
	
	\textbf{Appendix A} contains the proofs of the technical results used for the normal forms around homogeneous Lagrangian submanifolds: homogeneous Poincar\'e Lemma and tubular neighborhoods, and the proofs of the homogeneous Darboux-Weinstein, Weinstein Lagrangian neighborhood, and Weinstein tubular neighborhood Theorems. \textbf{Appendix B} derives a homogeneous version of the Hamilton-Jacobi equation.
	
	
	\subsubsection*{Acknowledgments} The authors would like to thank Oscar Cosserat for fruitful discussions and suggestions related with the work of this paper.
	The authors were financially supported by the Funda\c c\~ao para a Ci\^encia e a Tecnologia (FCT) under the scope of the project UID/00324 -- Center for Mathematics of the University of Coimbra. Gon\c calo Inoc\^encio Oliveira acknowledges FCT for support under the Ph.D.
	Scholarship 2024.00328.BD.
	
	
	\section{Homogeneous geometric structures}\label{anexo: homogeneous}

	As mentioned in the introduction, Jacobi manifolds can be viewed as homogeneous Poisson manifolds. In this section, we recall the setting for this interpretation: that of principal $\mathbb{R}^\times$-bundles and equivariant maps between them; we use the notation $\mathbb{R}^\times$ for the multiplicative group of non-zero reals.

	\begin{definition}[Principal  $\mathbb{R}^\times$-bundle]
		A principal $\mathbb{R}^\times$-bundle over a manifold $M$ consists of a manifold $P$ together with 
		\begin{itemize}
			\item[(i)] a right-action of $\mathbb{R}^\times$ on $P$, denoted by $h: P\times\mathbb{R}^\times\to P$, $(p,z)\mapsto h_z(p),$
			\item[(ii)] a surjective map $\tau:P\to M$ which is $\mathbb{R}^\times$-invariant, ($\tau(h_z(p))=\tau(p)$ for all $p$ and $z$),
		\end{itemize}
		satisfying local triviality: Every point $x_0\in M$ has an open neighborhood  $\mathcal{U}$ such that there is an $\mathbb{R}^\times$-equivariant diffeomorphism (called a local trivialization) $\psi_\mathcal{U}:\tau^{-1}(\mathcal{U})\rightarrow \mathcal{U}\times \mathbb{R}^\times$ which maps each fiber $\tau^{-1}(x)$ to the fiber $\{x\}\times \mathbb{R}^\times$. The action of $\mathbb{R}^\times$ on $\mathcal{U}\times \mathbb{R}^\times$ is by multiplication on the second factor. 
	\end{definition}
	
	We recall that the action of $\mathbb{R}^\times$ on $P$ is free and proper, and $\tau$ can be identified with the quotient map with respect to the action. 
	
	{
		\begin{remark}[On nomenclature]
			We will also refer to the pair $(P,h)$ as an \textit{$\mathbb{R}^\times$-manifold}. Although such objects are often called \textit{homogeneous} manifolds, we avoid this terminology to prevent confusion with homogeneous spaces; indeed, the $\mathbb{R}^\times$-action is not transitive except in trivial cases. Nevertheless, we retain the term \textit{homogeneous} when referring to related structures and objects compatible with the $\mathbb{R}^\times$-action, such as homogeneous maps, homogeneous symplectic manifolds, homogeneous Poisson manifolds, etc., in order to remain consistent with the terminology commonly used in the Jacobi geometry literature (cf.~\cite{VW2020}). In practice, the objects under consideration are simply principal bundles together with equivariant or invariant structures. However, in the specific case where the structure group is $\mathbb{R}^\times$, the \textit{homogenization scheme} (see Section~2 of \cite{VW2020}) provides a particularly powerful perspective on Jacobi geometry.
			
			The approach to homogenization adopted in this paper most closely follows that of \cite{BruceGrabowskaGrabowski2017}, which offers a streamlined treatment inspired by -- yet simpler than --  the graded manifold and graded bundle frameworks developed in \cite{GR2009,GR2012}. Earlier notions of homogeneity were sometimes formulated in terms of the associated infinitesimal action rather than the global $\mathbb{R}^\times$-action; such infinitesimal structures are inherited by open subsets. It is important to emphasize that, in our setting, as in \cite{BruceGrabowskaGrabowski2017,GR2012},  the infinitesimal generator of the action is required to be complete in order to integrate to an $\mathbb{R}^\times$-action. Consequently, when restricting to open subsets (for instance, coordinate chart domains), one must ensure that these subsets are saturated, i.e., invariant under the $\mathbb{R}^\times$-action.
		\end{remark}
	}

	\begin{definition}[Homogeneous map]
		Let $(P_1,h^1)$ and $(P_2,h^2)$ be $\mathbb{R}^\times$-manifolds. A smooth map $\phi:P_1\to P_2$ is called homogeneous if it is equivariant, i.e. if $\phi\circ h^1_z=h^2_z\circ \phi$, for all $z\in \mathbb{R}^\times$.
	\end{definition}
	
	\begin{example}[Frame bundle of a line bundle]
		Let $\tau_0:L\rightarrow M$ be a line bundle, i.e. a vector bundle of rank $1$.
		Consider the manifold $L^\times = L\backslash\{0_M\}$ which is the line bundle but without the zero section. We can define an action of multiplication of elements of $L^\times$ by non-zero reals as
		\begin{align*}
			h: L^\times \times \mathbb{R}^\times \rightarrow L^\times,\quad h(v,z)=h_z(v) = v\cdot z.
		\end{align*}
		Moreover, $\mathbb{R}^\times = GL_1(\mathbb{R})$ and each element $v\in L^\times_x$ forms a basis of the vector space $L_x$, for $x\in M$. So $L^\times = \text{Fr}(L)$ and $( L^\times,h,\tau)$ is a principal $\mathbb{R}^\times$-bundle, where $\tau : L^\times\rightarrow M$. If $(x,t)$ are coordinates in $L^\times$ given by a local trivialization of $L$, then the action $h$ is given by
		\[h_z(x,t) = (x,zt).\]
		In fact, up to isomorphism any $\mathbb{R}^\times$-bundle arises in this way from some associated line bundle.
	\end{example}

	Consider the principal $\mathbb{R}^\times$-bundle $\tau: P\rightarrow M$, where $P=L^\times$ as in the previous example, to simplify the notation. We can canonically lift the principal $\mathbb{R}^\times$-action on $P$ to principal $\mathbb{R}^\times$-actions on $TP$ (tangent lift) and $T^*P$ (phase lift). Consider the action defined previously $h_z(v) = v\cdot z$; then the tangent and phase lifts are
	\begin{align*}
		(Th)_z = Th_z\quad \text{ and }\quad  (T^*h)_z = z\cdot(Th_{z^{-1}})^*,
	\end{align*}
	see \cite{BruceGrabowskaGrabowski2017,Grabowska_2022}.
	
	In a local trivialization with coordinates $(x^i,t)$ on $P$, the $\mathbb{R}^\times$-action is $h_z(x^i,t) = ( x^i,zt)$. The natural local coordinates on $TP$ are $(x^i,t, \dot{x}^j, \dot{t})$, and so the tangent lift $Th$ is the action
	\begin{align}\label{eq: tangent lift}
		(Th)_z( x^i,t, \dot{x}^j, \dot{t}) = ( x^i,zt, \dot{x}^j, z\dot{t}).
	\end{align}
	We denote by $\dot{h}$ the action of $Th$ only on the tangent fibers, that is $\dot{h}_z(\dot{x},\dot{t}) = (\dot{x},z\dot{t}).$
	Similarly, the phase lift $T^*h$ acts on the cotangent coordinates $( x^i, t\xi_{x^j},\xi_t,)$ as
	\begin{align}\label{eq: cotangent lift}
		T^*h_z(x^i,t, \xi_{x^j},\xi_t) = (x^i,zt, z\xi_{x^j},\xi_t).
	\end{align}
	{\begin{remark}
			Hereinafter, whenever we take some coordinates $(x,t)$, $(x,t,\dot{x},\dot{t})$ or $(x^i,t, \xi_{x^j},\xi_t)$, consider always local coordinates saturated in $t$, $\dot{t}$ and $\xi_x$.
	\end{remark}}

	We can now define homogeneity for vector fields and differential forms. Denote by $\mathfrak{X}^m(M)$ the space of $m$-multivector fields on $M$, i.e., sections of $\Lambda^m(TM)$.
	
	\begin{definition}[$k$-homogeneous differential forms and multivector fields]
		Let $(P,h)$ be a $\mathbb{R}^\times$-manifold, let $\omega\in \Omega^l(P)$ and let $X\in \mathfrak{X}^m(P)$. We say that $\omega$ or $X$ are $k$-homogeneous (or homogeneous of degree $k$) if
		\begin{align}\label{form k-homoge}
			h_z^*\omega = z^k \omega
		\end{align}
		and respectively
		\begin{align}\label{vf k-homoge}
			(h_z)_*X = z^kX.
		\end{align}	
	\end{definition}
	
	We will also make use of submanifolds that respect the homogeneity structure.
	
	\begin{definition}[$\mathbb{R}^\times$-submanifolds]
		Let $(P,h)$ be a $\mathbb{R}^\times$-manifold and let $S$ be a submanifold of $P$. We say that $S$ is a $\mathbb{R}^\times$-submanifold (or homogeneous submanifold) if for every point $p\in S$, also $h_z(p)\in L$ holds.
	\end{definition}
	
	
	\section{Jacobi manifolds and Poissonization}\label{sec: poissonization}

	In this section, we give a short introduction to Jacobi manifolds and see how they can be interpreted as homogeneous Poisson manifolds by Poissonization. In the particular case of contact manifolds, seen as Jacobi manifolds, the Poissonization produces homogeneous symplectic manifolds.  A textbook account on Poisson geometry can be found in \cite{crainic2021lectures}, and on contact and Jacobi geometries in \cite{libermann1987symplectic}.
	\subsection{Poisson structures}\label{subsec: poisson-structures}
	
	\begin{definition}[Poisson structure]
		A Poisson structure on a differentiable manifold $M$ consists of a Lie bracket $\{\cdot,\cdot\}$ on the space $C^\infty (M)$ of smooth functions on $M$ satisfying additionally the Leibniz rule, i.e. it is a derivation in each entry \[\{f_1,f_2f_3\}=\{f_1,f_2\}f_3+f_2\{f_1,f_3\},\] for all $f_1,f_2,f_3\in C^\infty(M)$.
	\end{definition}
	
	A Poisson structure can equivalently be described as a bivector field $\Pi\in\mathfrak{X}^2(M)$ satisfying $[\Pi,\Pi]=0$, where $[\cdot,\cdot]$ is the Schouten-Nijenhuis bracket of multivector fields \cite[Chapter 2]{crainic2021lectures}.

	\begin{example}[Symplectic manifolds] Let $(M,\omega)$ be a symplectic manifold, meaning that $\omega\in\Omega^2(M)$ is closed and non-degenerate. 
		Then $M$ has a Poisson structure given by the canonical Poisson bracket from classical mechanics, $\{f_1,f_2\}:=\omega(X_{f_2},X_{f_1})$. Here $X_f$ denotes the Hamiltonian vector associated with $f$, uniquely defined by $\omega(X_f,\cdot)=df$. The associated Poisson bivector is denoted by $\Pi=\omega^{-1}$.
	\end{example}
	
	Given $H\in C^\infty(M)$, the operation $\{H,\cdot\}$ is a derivation of $C^\infty(M)$, so it is a vector field, denoted by $X_H$ and called the \textit{Hamiltonian vector field} associated with $H$. The vector subspaces of the tangent spaces to $M$ given by the value of all possible vector fields form a smooth distribution on $M$. Although this distribution is singular in general, in the sense that it might have different dimensions at different points, it is integrable: there is a partition of $M$ into submanifolds forming a singular foliation. These submanifolds, called \textit{symplectic leaves}, carry symplectic structures induced by the Poisson structure.
	
	This symplectic foliation is relevant to understand the qualitative aspects of Hamiltonian dynamics on a Poisson manifold. The symplectic leaf containing a point $p$ is composed of all the points that can be reached by starting from $p$, by repeatedly following the flow of Hamiltonian vector fields. 
	For example, the symplectic leaves of a symplectic manifold $(M,\omega)$ seen as a Poisson manifold are just the connected components of $M$.
	
	
	\subsection{Jacobi structures and examples}
	\begin{definition}[Jacobi structure \cite{LICHNEROWICZ1977}]
		Let $J$ be a smooth manifold and let $\Lambda$ be a bivector field and $E$ a vector field on $J$, respectively. We call $(\Lambda, E)$ a Jacobi structure if
		\begin{align}\label{eq: jacobi}
			[\Lambda,\Lambda]= 2E\wedge\Lambda,\qquad [\Lambda,E] = 0,
		\end{align}
		where $[\cdot,\cdot]$ is the Schouten-Nijenhuis bracket. We call the triple $(J,\Lambda,E)$ a Jacobi manifold.
	\end{definition}
	
	Associated with a Jacobi manifold $(J,\Lambda,E)$, we can define a Jacobi bracket by
	\begin{align}\label{eq: J bracket}
		\{f_1,f_2\}_J = \Lambda(df_1,df_2) + f_1E(f_2) - f_2E(f_1),\quad f_1,f_2\in C^\infty(J).
	\end{align}
	
	This is a Lie bracket on the space of smooth functions on $J$ which satisfies
	\begin{align*}
		\{f_1f_2,f_3\}_J = f_1\{f_2,f_3\}_J +f_2\{f_1,f_3\}_J - f_1f_2\{1,f_3\}_J,\quad f_1,f_2,f_3\in C^\infty(J).
	\end{align*}
	
	Note that, from the last two terms we have a vector field associated with $H\in C^\infty(J)$,
	\begin{align*}
		X_H = \Lambda(\cdot,dH) - HE(\cdot),
	\end{align*}
	which is called the Hamiltonian vector field associated with $H$. For example, $E$ is the Hamiltonian vector field associated with the constant function $-1$.

	\begin{remark}[On the general definition of Jacobi manifold]\label{rmk-Jacobi}
		The definition of a Jacobi structure adopted here is not the most general one. In full generality, Jacobi brackets are defined on the module of sections of a line bundle, following the approach of Kirillov (cf.~\cite{Kirillov,BruceGrabowskaGrabowski2017}). For simplicity of exposition, we restrict ourselves to the case of a trivial line bundle. This setting is sufficient for the purposes of the present work, since the construction and implementation of geometric numerical integrators for Hamiltonian systems can be carried out within a local trivialization.
	\end{remark}
	
	Jacobi structures can be seen as common generalizations of contact structures as well as of Poisson structures.
	
	\begin{example}[Poisson]
		Let $(M,\Pi)$ be a Poisson manifold. We can interpret it as a Jacobi manifold letting $\Lambda = \Pi$ and $E =0$. In this case, Equation (\ref{eq: jacobi}) amounts to the definition of a Poisson structure, $[\Pi,\Pi] = 0$.
	\end{example} 
	
	Another important example of Jacobi manifolds are contact manifolds, which we now present following \cite{da2001lectures}.
	
	\begin{definition}[Contact manifold]
		Let $M$ be a $(2n+1)$-smooth manifold. A contact structure on $M$ is a distribution of hyperplanes $\mathfrak{H} \subset TM$, maximally non-integrable, for which there exists locally a $1$-form $\eta$ such that $\mathfrak{H} = \text{ker}\ \eta$ and $d\eta_{|\mathfrak{H}}$ is nondegenerate (i.e., symplectic). The pair $(M, \mathfrak{H})$ is then called a contact manifold and $\eta$ is called a local contact form. If $\mathfrak{H} = \text{ker}\ \eta$ globally, we call $\eta$ a contact form.
	\end{definition}
	
	Given a contact form $\eta$, there is an associated Reeb vector field $\xi$, defined as the unique vector field satisfying
	\[
	\eta(\xi) = 1 \qquad \text{and} \qquad i_\xi d\eta = 0.
	\]
	Furthermore, $\eta$ induces an isomorphism
	\begin{align*}
		\flat_\eta: TM &\longrightarrow T^*M,\\
		X &\longmapsto i_X\eta \,\eta + i_X d\eta.
	\end{align*}

	\begin{example}[Contact \cite{Simoes_2020}]\label{ex: contact}
		Let $(M,\mathfrak{H})$ be a contact manifold with local contact form $\eta$ and let $\xi$ be the associated Reeb vector field. Using the isomorphism $\flat_\eta$, we define the bivector field as $\Lambda(\alpha,\beta) = -d\eta \left(\flat_\eta^{-1}(\alpha), \flat_\eta^{-1}(\beta)\right)$ and the vector field $E=-\xi$. In canonical coordinates it takes the form \begin{align*}
			\Lambda = \frac{\partial}{\partial p_i}\wedge\frac{\partial}{\partial q^i} + p_i\frac{\partial}{\partial p_i}\wedge\frac{\partial}{\partial z}.
		\end{align*}
		and $E=-\frac{\partial}{\partial z}$. 
		When the contact form is global, this defines a Jacobi structure on $M$. To allow for a similar construction to define a Jacobi structure in the general case, Jacobi structures supported in non-trivial line bundles (as mentioned in Remark \ref{rmk-Jacobi}) are needed.
	\end{example}
	
	The main contribution of this work is a systematic method to construct structure-preserving integrators for Hamiltonian vector fields on Jacobi manifolds. A key initial step in our approach is the passage from Jacobi manifolds to Poisson manifolds through Poissonization.
	
	
	\subsection{Poissonization and homogeneous Poisson manifolds}
	
	The process of Poissonization from \cite{LICHNEROWICZ1977} translates Jacobi manifolds into homogeneous Poisson manifolds, thereby embedding Jacobi geometry into the broader Poisson framework. This construction is both explicit and canonical, and it plays a central role in understanding Jacobi structures as Poisson structures with an additional scaling symmetry.
	
	\begin{definition}[Poissonization]
		Let $(J, \Lambda, E)$ be a Jacobi manifold. Let $P_J := J \times \mathbb{R}^\times$ be the $\mathbb{R}^\times$-principal bundle over $J$ given by the principal $\mathbb{R}^\times$-action $h_z(x,t) = (x, zt)$. Consider the bivector field $\Pi$ on $P_J$ given by
		\[
		\Pi(x, t) = \frac{1}{t} \Lambda(x) + \frac{\partial}{\partial t} \wedge E(x),
		\]
		along with the vector field on $P_J$
		\[
		Z := t \frac{\partial}{\partial t},
		\]
		which is the infinitesimal generator of the principal $\mathbb{R}^\times$-action. The triple $(P_J, \Pi, h)$ is called the Poissonization of $(J, \Lambda, E)$.
	\end{definition}
	
	We now formalize the notion of a homogeneous Poisson or symplectic manifold.
	
	\begin{definition}[Homogeneous Poisson and symplectic manifolds]
		Let $(P, \pi)$ be a Poisson manifold equipped with a principal $\mathbb{R}^\times$-action $h_z$. We say that $(P, \pi, h)$ is a homogeneous Poisson manifold if $\pi$ is $-1$-homogeneous, i.e.
		\[
		(h_z)_* \pi = \frac{1}{z} \pi \quad \text{for all } z \in \mathbb{R}^\times.
		\]
		Similarly, a symplectic manifold $(\Sigma, \omega)$ with an action $h_z$ is called a \emph{homogeneous symplectic manifold} if $\omega$ is $1$-homogeneous, i.e.
		\[
		h_z^* \omega = z \omega \quad \text{for all } z \in \mathbb{R}^\times.
		\]
	\end{definition}
	
	\begin{proposition}
		The Poissonization $(P_J, \Pi, h)$ defined above is a homogeneous Poisson manifold with respect to the principal action $h_z(x,t) = (x, zt)$.
	\end{proposition}
	
	\noindent \textit{Proof:}
		A direct computation shows that
		\[
		(h_z)_* \Pi = \frac{1}{z} \Pi,
		\]
		so the defining condition for homogeneity holds.  
		
		\noindent $\square$
	
	\begin{remark}
		In the case where $(J, \Lambda, E)$ is a contact manifold, the associated Poissonization yields a homogeneous symplectic manifold, known as its \emph{symplectization}. The Poisson bivector $\Pi$ is then non-degenerate, and the symplectic form $\omega=\Pi^{-1}$ satisfies $h_z^* \omega = z\omega$.  Thus, the symplectic category can be seen as a homogeneous lift of contact geometry.
	\end{remark}
	
	The Poissonization process reflects a deep equivalence between categories.
	
	\begin{proposition}[\cite{Monica_Blaga_2020}, Proposition B.5]\label{Prop-equiv-cat}
		There is an equivalence of categories between the Jacobi (resp. contact) category and the homogeneous Poisson (resp. symplectic) category.
	\end{proposition}
	
	\subsection{Casimir functions}
	
	As an illustration of the use of Poissonization, in this subsection we study Casimir functions, invariants that are annihilated by the bracket structure of a manifold. In the Jacobi setting, these functions are characterized by vanishing both under the image of the Jacobi bivector $\Lambda$ and under the vector field $E$. We examine how these functions lift naturally through the Poissonization process, becoming homogeneous Poisson-Casimirs.

	\begin{definition}[Jacobi-Casimir functions]
		Let $(J, \Lambda, E)$ be a Jacobi manifold. A smooth function $f\in C^\infty(J)$ is called a Jacobi-Casimir function if
		\[
		\Lambda^\sharp(df) = 0 \quad \text{and} \quad E(f) = 0.
		\]
	\end{definition}
	
	\begin{remark}
		When $E = 0$, the condition reduces to the usual notion of a Poisson-Casimir: $\Pi^\sharp(df) = 0$.
	\end{remark}
	
	Let $f \in C^\infty(P_J)$ be a $0$-homogeneous function, that is, $f(x,t) = f(x)$ does not depend on $t$. Then
	\begin{align*}
		\Pi^\sharp(df) &= \left( \frac{1}{t} \Lambda + \frac{\partial}{\partial t} \wedge E \right)^\sharp(df) \\
		&= \frac{1}{t} \Lambda^\sharp(df) + df(E) \frac{\partial}{\partial t} - df\left(\frac{\partial}{\partial t}\right) E \\
		&= \frac{1}{t} \Lambda^\sharp(df) + E(f) \frac{\partial}{\partial t},
	\end{align*}
	since $df$ is independent of $t$. Thus, $\Pi^\sharp(df) = 0$ if and only if $f$ is a Jacobi-Casimir function.
	
	\begin{remark}
		We conclude that Jacobi-Casimir functions lift to $0$-homogeneous Poisson-Casimir functions on the Poissonized space $(P_J, \Pi)$. This compatibility confirms the coherence of the Poissonization construction with the underlying algebraic structures.
	\end{remark}
	

	\section{Tools from homogeneous symplectic and Poisson geometry}\label{sec-homogeneous-tools}
	
	This section describes the geometric tools which will be used in our construction of structure-preserving integrators in the Jacobi setting, namely:
	\begin{enumerate}\item results on Lagrangian and Legendrian manifolds (key ingredients in the geometric Hamilton-Jacobi theory);
		\item explicit constructions of (homogeneous) symplectic realizations.
	\end{enumerate}
	
	
	\subsection{Legendrian and homogeneous Lagrangian submanifolds}
	
	Through the interpretation of contact manifolds and homogeneous symplectic manifolds, we now consider some results about special submanifolds in contact geometry. The first result relates homogeneous Lagrangian submanifolds \cite{vogtmann1997mathematical} and Legendrian submanifolds \cite{Eliashberg1989}.
	
	\begin{definition}[Homogeneous Lagrangian submanifold]
		Let $(\Sigma,\omega)$ be a $2n$-dimensional symplectic manifold and let $L$ be an $n$-dimensional submanifold. We say that $L$ is a Lagrangian manifold if $i^*\omega=0$, where $i:L\xhookrightarrow{} \Sigma$ is the inclusion map. That is, the symplectic form vanishes on vectors tangent to $L$. If additionally $(\Sigma,\omega,h)$ is a homogeneous symplectic manifold, L is called a homogeneous Lagrangian submanifold if it is both homeogeneous (i.e. $\mathbb{R}^\times$-invariant) and Lagrangian.
	\end{definition}
	
	In contact geometry, there is a notion similar to  Lagrangian submanifolds, that of Legendrian submanifolds.
	
	\begin{definition}[Legendrian submanifold]
		Let $(C,\mathfrak{H})$ be a $(2n+1)$-dimensional contact manifold with contact form $\eta$ and let $L$ be an $n$-dimensional submanifold. We call $L$ a Legendrian submanifold if for every $p\in C$, $T_pL\in \mathfrak{H}_p$, that is $T_pL\in \textnormal{ker}\ \eta_p$. 
	\end{definition}

	Now, consider a contact manifold $(C,\eta)$ and its symplectization $(\Sigma,\omega,h)$.
	
	\begin{proposition}[\cite{Grabowska_2022}]\label{prop: legendrian and lagrangian} 
		There is a canonical one-to-one correspondence between $\mathbb{R}^\times$-invariant (or homogeneous) Lagrangian submanifolds $\mathcal{L}$ of $\Sigma$ and Legendre submanifolds $\mathcal{L}_0=\tau(\mathcal{L})$ of $C$.
	\end{proposition}
	
	Lagrangian submanifolds are of very wide utility in symplectic geometry; we will make use of the homogeneous version of them to codify both 1-forms and maps, through the following two results. The first one relates the image of homogeneous closed 1-forms to homogeneous Lagrangian submanifolds. Let $\mu$ be a 1-homogeneous 1-form on $\Sigma$, a homogeneous symplectic manifold, and consider its image $X_\mu = \{(x,\mu_x)|x\in \Sigma, \mu_x\in T^*_x\Sigma\}$.
	
	\begin{proposition}\label{prop: mu closed}
		$X_\mu$ is a homogeneous Lagrangian of $T^*\Sigma$ if and only if $\mu$ is a 1-homogeneous closed $1$-form.
	\end{proposition}
	
	\noindent \textit{Proof:}
		From \cite[Chapter~3]{da2001lectures} we know that $X_\mu$ is Lagrangian if and only if $d\mu =0$. We are only left with proving the homogeneity.
		
		Suppose that {$(x,t)$ are homogeneous coordinates} of $\Sigma$. We know that $\mu$ is $1$-homogeneous if $h_z^*\mu = z\mu$. Consider $\mu = f(x,t)dx + g(x,t)dt$. In this expression, the homogeneity of $\mu$ is equivalent to $h_z^*f(x,t) = zf(x,t)$ and $h_z^*g(x,t) = g(x,t)$. 
		
		Taking into account the embedding $s_\mu:\Sigma\rightarrow T^*\Sigma,\ x\mapsto (x,\mu_x)$, its image is precisely $X_\mu$. So, the homogeneity of $X_\mu$ is related to the homogeneity of $s_\mu$. In the previous identification, $s_\mu(x,t) = (x,t,f(x,t),g(x,t))$, and with the lifted action (\ref{eq: cotangent lift}) we get
		\begin{align*}
			(s_\mu\circ h_z)(x,t) & = (x,zt,f(x,zt),g(x,zt))\\
			& = (x,zt,zf(x,t),g(x,t))\\
			& = (T^*h_z)(x,t,f(x,t),g(x,t))\\
			& = (T^*h_z\circ s_\mu)(x,t). 
		\end{align*}
		
	\noindent $\square$ 
	
	We now  consider two homogeneous symplectic manifolds $(\Sigma_1,\omega_1,h^1)$ and $(\Sigma_2, \omega_2, h^2)$, so that $(h_z^i)^*\omega_i = z\omega_i,\ i=1,2$. Given a homogeneous diffeomorphism $\varphi:\Sigma_1\rightarrow \Sigma_2$, the next proposition characterizes when it is a homogeneous symplectomorphism.
	
	\begin{proposition}\label{prop: phi lagrangian}
		A homogeneous diffeomorphism $\varphi: \Sigma_1\rightarrow \Sigma_2$ is a homogeneous symplectomorphism if and only if its graph $\Gamma_\varphi$ is a homogeneous Lagrangian submanifold of $(\Sigma_1\times \Sigma_2,\pi_1^*\omega_1 - \pi_2^*\omega_2, h^\times)$.
	\end{proposition}
	
	\noindent \textit{Proof:}
		Consider $\Sigma_1\times \Sigma_2$ with projections $\pi_i:\Sigma_1\times \Sigma_2 \rightarrow \Sigma_i,\ i=1,2$ and also the symplectic form $\omega = \pi_1^*\omega_1 - \pi_2^*\omega_2$. We know by \cite[Proposition 3.8]{da2001lectures} that $\varphi$ is a symplectomorphism if and only if $\Gamma_\varphi$ is Lagrangian.
		
		Consider the multiplication $h^\times_z = (h^1\times h^2)_z$ on $\Sigma_1\times \Sigma_2$. First, the symplectic form $\omega$ is 1-homogeneous
		\begin{align*}
			(h^\times_z)^*\omega &= (h^\times_z)^*(\pi_1^*\omega_1 - \pi_2^*\omega_2)\\
			& = (\pi_1\circ h^\times_z)^*\omega_1 - (\pi_2\circ h^\times_z)^*\omega_2\\
			& = (h^1_z\circ\pi_1)^*\omega_1 - (h^2_z\circ\pi_2)^*\omega_2\\
			& = \pi_1^*(h^1_z)^*\omega_1 - \pi_2^*(h^2_z)^*\omega_2\\
			& = z\pi_1^*\omega_1 - z\pi_2^*\omega_2\\
			& = z\omega.
		\end{align*}
		
		We also know that the graph $\Gamma_\varphi$ is an embedded image of $\Sigma_1$ in $\Sigma_1\times \Sigma_2$ with embedding $\gamma:\Sigma_1\rightarrow M_1\times M_2,\ p\mapsto (p,\varphi(p))$. This embedding is homogeneous if and only if $\varphi$ is homogeneous:
		\begin{align*}
			(h^\times_z\circ \gamma) (p) = h^\times_z(p,\varphi(p)) = \left(h^1_z(p), h^2_z(\varphi(p))\right) = \left(h^1_z(p), \varphi(h^1_z(p))\right) = (\gamma\circ h^1_z)(p).
		\end{align*}
		
	\noindent $\square$  
	
	Finally, we will make use of the two following theorems, in order to obtain local normal forms in homogeneous neighborhoods of submanifolds of a homogeneous symplectic manifold. 
	
	Normal forms are essentially good choices of coordinates in which the structures at hand can be described in a simpler way. This is reminiscent, for example, of the Jordan normal form in Linear Algebra, which is given by a choice of basis for which a linear map is described in a simple fashion.
	
	\begin{thm}[homogeneous Darboux-Weinstein theorem]\label{thm: Darboux-Weinstein}
		Let $(M,h)$ be a $\mathbb{R}^\times$-manifold and let $X\subset M$ be a $\mathbb{R}^\times$-submanifold. Suppose $\omega_0,\ \omega_1$ are two $1$-homogeneous symplectic forms on $M$, for which $(\omega_0)_{|X} =(\omega_1)_{|X} $. Then, there is a { homogeneous neighborhood $\mathcal{U}$} of $X$ and a diffeomorphism $f:\mathcal{U}\rightarrow \mathcal{U}$ such that:
		\begin{enumerate}
			\item $f(x) = x$, for all $x\in X$;
			\item $f^*\omega_1 = \omega_0$;
			\item if $h_z$ is the principal action, then $f\circ h = h\circ f.$
		\end{enumerate}
	\end{thm}

	\begin{thm}[Homogeneous Weinstein tubular neighborhood theorem]\label{thm-WLtub} 	
		Let $(M,\omega, h)$ be a $2n$-dimensional homogeneous symplectic manifold, let $X$ be an homogeneous Lagrangian submanifold, with $i:X\xhookrightarrow{} M$ the inclusion map, and $i_0:X\xhookrightarrow{} T^*X$ the Lagrangian embedding as the zero section. Let $\omega_0$ be the canonical symplectic form on $T^*X$. Then, there exist homogeneous neighborhoods $\mathcal{U}_0$ of $X$ in $T^*X$ and $\mathcal{U}$  of $X$ in $M$ and a homogeneous diffeomorphism $\varphi:\mathcal{U}_0\rightarrow \mathcal{U}$ such that $\varphi^*\omega = \omega_0$ and the following diagram commutes
		\begin{center}
			\begin{tikzpicture}
				\node at (0,0) (a) {$\mathcal{U}_0$};
				\node at (1,-1) (b) {$X$};
				\node at (2,0) (c) {$\mathcal{U}$};
				\draw [->] (b) to node[left] {$i_0$} (a);
				\draw [->] (a) to node[above] {$\varphi$} (c);
				\draw [->] (b) to node[right] {$i$}(c);
			\end{tikzpicture}
		\end{center}
	\end{thm}
	
	We give the proofs of these two theorems in Appendix \ref{anexo: symplectic}. The key steps, included in the Appendix, are to use homogeneous versions of the Poincar\'e Lemma, of tubular neighborhoods, of the Moser trick, and of the Weinstein Lagrangian neighborhood Theorem.
	
	
	\subsection{Constructing homogeneous bi-realizations}\label{sec: bi-realization}
	
	We will now define, and explicitly construct homogeneous symplectic bi-realizations. In Poisson geometry, one way to construct symplectic realizations is through an auxilliary Poisson spray. We prove that a Poisson spray can be chosen such that it preserves the homogeneity property, which leads to homogeneous bi-realizations.
	
	\begin{definition}[Symplectic realization] A symplectic realization of a Poisson manifold $(M,\pi)$, denoted by
		\begin{align*}
			\mu:(S,\omega) \rightarrow(M,\pi),
		\end{align*}
		consists of:
		\begin{enumerate}
			\item a symplectic manifold $(S,\omega)$;
			\item a surjective submersion $\mu:S\rightarrow M$ which is a Poisson map.
		\end{enumerate}
	\end{definition}
	
	\begin{definition}[Bi-realization] Let $\pi$ be a Poisson structure on an open subset $U\subset \mathbb{R}^n$. A bi-realization of $(U,\pi)$ is given by a bi-surjection, denoted by $W\rightrightarrows U$, \i.e., a pair of surjective submersions $\alpha$ and $\beta$, called source and target, respectively, satisfying:
		\begin{enumerate}
			\item $\alpha$ is a Poisson map;
			\item $\beta$ is an anti-Poisson map;
			\item the fibers of $\alpha$ and $\beta$ are symplectic orthogonal to each other.
		\end{enumerate}	
	\end{definition}

	By \cite{Monica_Blaga_2020} we have the following definition of homogeneous symplectic bi-realization.
	
	\begin{definition}[Homogeneous symplectic bi-realization]
		A homogeneous symplectic bi-realization is a symplectic bi-realization $(\Sigma,\omega)\rightrightarrows (P,\{\cdot,\cdot\})$ such that the manifolds $\Sigma$ and $P$ are equipped with principal $\mathbb{R}^\times$-bundles structures, and both the Poisson and the symplectic structures, as well as source and target maps, are homogeneous. 
	\end{definition}
	
	One way to construct explicit bi-realizations for Poisson manifolds is using a Poisson spray \cite{crainic2021lectures}. For a homogeneous Poisson manifold, what we need to prove is that the same construction is compatible with homogeneity, identifying the correct choices of spray for that to happen.

	\begin{definition}[Poisson spray]
		Let $P$ be a Poisson manifold. A Poisson spray is a vector field $X\in\mathfrak{X}(T^*P)$ that satisfies the following:
		\begin{enumerate}
			\item[(i)] $d_\xi\ \tau(X_\xi)= \pi^\sharp(\xi)$, for all $\xi\in T^*P$,
			\item[(ii)] $(m_t)_*X = \frac{1}{t}X$, for all $t>0$,  
		\end{enumerate}
		where $\tau:T^*P\rightarrow P$ denotes the cotangent projection and $m_t:T^*P\rightarrow T^*P$ is the scalar multiplication by $t\in \mathbb{R}$.
	\end{definition}
	
	Let $y$ be coordinates on $P$. Then, the Poisson spray can be written locally as
	\begin{align*}
		X &= \sum_{ij}\Pi_{ij}(y)\xi_j\frac{\partial}{\partial y_i} + f_{ij}\xi_i\xi_j\frac{\partial}{\partial\xi_j}.
	\end{align*}
	
	In our case, if $(x,t)$ are the coordinates of $P_J$, the Poisson spray takes the form
	\begin{align}\label{eq: spray}
		X = \sum_{ij}\frac{1}{t}\Lambda_{ij}(x)\xi_j\frac{\partial}{\partial x_i} + E_j(x)\xi_j\frac{\partial}{\partial t} - E_i(x)\xi_t\frac{\partial}{\partial x_i} + f_{ij}\xi_i\xi_j\frac{\partial}{\partial\xi_j},
	\end{align}
	where $f_{ij}\in C^\infty(P),\ i,j\in\{1,\dots,2n+1,t\}$ are functions free to choose.
	
	If we compute the homogeneity of the non-optional terms of $X$ we see that
	\begin{align*}
		(T^*h_z)_*X &= (T^*h_z)_*\left(\sum_{ij}\frac{1}{t}\Lambda_{ij}(x)\xi_j\frac{\partial}{\partial x_i} + E_j(x)\xi_j\frac{\partial}{\partial t} - E_i(x)\xi_t\frac{\partial}{\partial x_i} \right)\\
		&=\sum_{ij}\frac{1}{zt}\Lambda_{ij}(x)z\xi_j\frac{\partial}{\partial x_i} + E_j(x)z\xi_j\frac{1}{z}\frac{\partial}{\partial t} - E_i(x)\xi_t\frac{\partial}{\partial x_i} \\
		&=\sum_{ij}\frac{1}{t}\Lambda_{ij}(x)\xi_j\frac{\partial}{\partial x_i} + E_j(x)\xi_j\frac{\partial}{\partial t} - E_i(x)\xi_t\frac{\partial}{\partial x_i} 
	\end{align*}
	all terms are $0$-homogeneous. We want $X$ to be $0$-homogeneous, so we need to find the conditions for the optional terms
	\begin{align*}
		&(T^*h_z)_*\left(f_{ij}(x,t)\xi_i\xi_j\frac{\partial}{\partial\xi_j} + f_{tj}(x,t)\xi_t\xi_j\frac{\partial}{\partial\xi_j}+f_{it}(x,t)\xi_i\xi_t\frac{\partial}{\partial\xi_t} + f_{tt}(x,t)\xi_t\xi_t\frac{\partial}{\partial\xi_t}\right) \\
		=& f_{ij}(x,zt)z^2\xi_i\xi_j\frac{1}{z}\frac{\partial}{\partial\xi_j} + f_{tj}(x,zt)\xi_tz\xi_j\frac{1}{z}\frac{\partial}{\partial\xi_j}+f_{it}(x,zt)z\xi_i\xi_t\frac{\partial}{\partial\xi_t} + f_{tt}(x,zt)\xi_t\xi_t\frac{\partial}{\partial\xi_t}\\
		=& f_{ij}(x,zt)z\xi_i\xi_j\frac{\partial}{\partial\xi_j} + f_{tj}(x,zt)\xi_t\xi_j\frac{\partial}{\partial\xi_j}+f_{it}(x,zt)z\xi_i\xi_t\frac{\partial}{\partial\xi_t} + f_{tt}(x,zt)\xi_t\xi_t\frac{\partial}{\partial\xi_t}.
	\end{align*}
	
	We can divide into two cases:
	\begin{align}\label{eq: f_ij}
		\begin{cases}
			f_{ij}(x,zt) = \frac{1}{z}f_{ij}(x,t),\ \textnormal{when } i\neq t,\\
			f_{tj}(x,zt) = f_{t,j(x,t)}.
		\end{cases} 
	\end{align}
	
	This means that, when $i\neq t$, the $f_{i,j}$ need to be $-1$-homogeneous and when $i=t$, the $f_{tj}$ need to be $0$-homogeneous.
	
	If (\ref{eq: f_ij}) holds, then $X$ is $T^*h_z$-related with itself, so 
	\[\phi_X^s\circ T^*h_z = T^*h_z\circ \phi_X^s,\]
	where $\phi_X$ is the flow of $X$.
	
	\begin{proposition}
		If (\ref{eq: f_ij}) holds, then the Poisson spray, given by (\ref{eq: spray}), is $0$-homogeneous and its flow commutes with $T^*h_z$.
	\end{proposition}
	
	Now we will construct explicit homogeneous bi-realizations in canonical form, (essentially similar to the construction of local integrations of Jacobi structures of \cite{CMS2022}). The technique of using a Poisson spray was developed in \cite{CrainicMarcut2011}, and in \cite{COSSERAT2023} the authors make use of the construction of bi-realizations using the same technique by \cite{CMS2022}.

	Consider the cotangent projection $\tau:T^*P_J\rightarrow P_J$. It is homogeneous $h_z\circ \tau = \tau \circ T^*h_z$, so
	\begin{align*}
		&\bar{\alpha} := \tau \longrightarrow \textnormal{ homogeneous}\\
		&\bar{\beta} := \tau\circ \phi^1_X \longrightarrow \textnormal{homogeneous}\\
		&\qquad \tau\circ\phi^1_X \circ T^*h_z = \tau\circ T^*h_z\circ \phi^1_X = h_z\circ\tau\circ\phi^1_X
	\end{align*}
	and since $\omega_{can} = dx\wedge d\xi_x + dt\wedge d\xi_t$ is $1$-homogeneous, the symplectic form $\Omega = \int_0^1(\phi^s_X)^*\omega_{can} ds$ is also $1$-homogeneous:
	\begin{align*}
		(T^*h_z)^*\Omega &= (T^*h_z)^*\int_0^1(\phi^s_X)^*\omega_{can}ds = \int_0^1(T^*h_z)^*(\phi^s_X)^*\omega_{can}ds\\
		&=\int_0^1(\phi^s_X\circ T^*h_z)^*\omega_{can}ds = \int_0^1 (T^*h_z\circ\phi^s_X)^*\omega_{can}ds\\
		& = \int_0^1 (\phi^s_X)^*(T^*h_z)^*\omega_{can}ds = z\int_0^1(\phi^s_X)^*\omega_{can}ds.
	\end{align*}

	If we stop here, we already have a homogeneous symplectic bi-realization. In this case, the realization maps are the cotangent projection and its composition with the flow of the spray, and the symplectic form is a deformation of the canonical symplectic form; this is called the Weinstein realization by \cite{cabrera2024}. However, for our construction of Jacobi Hamiltonian integrators it will be convenient to have instead a realization where the symplectic form is the canonical one and the realization maps are deformations of the cotangent projection, which is the Karasev realization \cite{Karasev_1987}.

	By the homogeneous Darboux-Weinstein Theorem (\ref{thm: Darboux-Weinstein}), there exists $\Psi:T^*P\rightarrow T^*P$ such that $\omega_{can} = \Psi^*\Omega$ and $(T^*h_z)^*\Psi = \Psi$ in a { homogeneous }neighborhood of the zero section in $T^*P$. So, we can define $\alpha = \Psi^*\bar{\alpha}$ and $\beta=\Psi^*\bar{\beta}$, obtaining a new bi-realization.
	
	\begin{proposition}\label{thm: bi-realization}
		Any homogeneous Poisson spray induces a homogeneous bi-realization.
	\end{proposition}
	
	\noindent \textit{Proof:}
		We know that any Poisson spray induces a bi-realization with source map $\alpha$ and target $\beta$. We only need to prove that these maps are homogeneous. From the homogeneity of $\Psi$ and $\phi^s_X$ we get
		\begin{align*}
			h_z\circ\alpha &= h_z\circ \tau \circ \Psi = \tau\circ T^*h_z \circ \Psi = \tau\circ\Psi\circ T^*h_z = \alpha\circ T^*h_z
		\end{align*}
		and
		\begin{align*}
			h_z\circ\beta = h_z\circ \tau\circ \phi^1_X\circ \Psi = \tau\circ\phi^1_X \circ T^*h_z\circ\Psi =\tau\circ\phi^1_X \circ\Psi\circ T^*h_z = \beta\circ T^*h_z.
		\end{align*}
		So, by definition, the bi-realization is homogeneous. 
	
	\noindent $\square$  
	
	\begin{remark}
		The Karasev realization has the property
		\[\beta(x,\xi)=\alpha(x,-\xi).\]
	\end{remark}

	As expected, there is a relation between homogeneous symplectic bi-realizations and contact realizations.
	
	\begin{proposition}[\cite{CrainicZhu2004,Monica_Blaga_2020}]\label{prop: hsb and cb}
		There exists a $1$-to-$1$ correspondence between homogeneous symplectic bi-realizations and contact bi-realizations.
	\end{proposition}
	
	
	\subsection{An example of a homogeneous bi-realization}\label{sec: example}
	
	In this section, we construct the Poissonization and the explicit bi-realization for the canonical contact case. As in Example \ref{ex: contact}, the canonical contact structure is given by
	\begin{align}\label{eq: contact structure example}
		\Lambda = \frac{\partial}{\partial p}\wedge\frac{\partial}{\partial q} + p\frac{\partial}{\partial p_i}\wedge\frac{\partial}{\partial z}\qquad \text{ and } \qquad E=-\frac{\partial}{\partial z}.
	\end{align}
	Doing the Poissonization trick, we have the following 
	\begin{align*}
		\Pi = \frac{1}{t}\frac{\partial}{\partial p}\wedge\frac{\partial}{\partial q} + \left(\frac{p}{t} \frac{\partial}{\partial p} - \frac{\partial}{\partial t} \right)\wedge \frac{\partial}{\partial z}
	\end{align*}
	and the flat Poisson spray is
	\begin{align*}
		X = \frac{1}{t}\left(\xi_{q} + p\xi_z\right)\frac{\partial}{\partial p} - \frac{1}{t}\xi_{p}\frac{\partial}{\partial q} + \left(\xi_t-\frac{p}{t}\xi_{p} \right)\frac{\partial}{\partial z} - \xi_z\frac{\partial}{\partial t}.
	\end{align*}
	However, if we compute its flow $\phi^s_X$, the component of $p$ is the expression
	\begin{align*}
		\phi^s_{p}(x,\xi) = pe^{\int_1^s\frac{\xi_z(\zeta)}{t(\zeta)}d\zeta} + e^{\int_1^s\frac{\xi_z(\zeta)}{t(\zeta)}d\zeta}\int_1^s e^{-\int_1^\zeta\frac{\xi_z(\delta)}{t(\delta)}d\delta}\frac{\xi_{q}(\zeta)}{t(\zeta)}d\zeta,
	\end{align*}
	which is impractical to work with. 
	
	{An alternative approach is to construct a bi-realization via a symplectomorphism. Consider
		\[
		F(q,p,z,t) = (Q,P,Z,T) = (q, pt, -z, t),
		\]
		which transforms the original Poisson structure $\Pi$ into the canonical symplectic structure
		\[
		F_*\Pi = \frac{\partial}{\partial P}\wedge \frac{\partial}{\partial Q} + \frac{\partial}{\partial T}\wedge \frac{\partial}{\partial Z} = \omega^{-1}_{\text{can}}.
		\]
		The canonical Poisson spray is
		\[X_\text{can}= \xi_Q\dfrac{\partial}{\partial P} - \xi_P\dfrac{\partial}{\partial Q} + \xi_Z\dfrac{\partial}{\partial T} - \xi_T\dfrac{\partial}{\partial Z},\]
		whose flow is given explicitly by
		\[\phi^s_{X_\text{can}} = (Q-s\xi_p, P+s\xi_Q, Z-s\xi_T, T+s\xi_Z, \xi_Q,\xi_P,\xi_Z,\xi_T).\]
		This flow is linear in the fiber coordinates, making computations straightforward. The associated canonical form is
		\[\Omega_\text{can} = \int_0^1 (\phi^s_{X_\text{can}})^*\omega_\text{can} ds = \omega_\text{can} +\frac{1}{2}(dQ\wedge d\xi_Q -d\xi_P\wedge dP + dZ\wedge d\xi_Z - d\xi_T\wedge dT) + \frac{1}{3}(d\xi_P\wedge \xi_Q - d\xi_T\wedge d\xi_Z).\]
		Using the symplectomorphism $\Psi = \phi^{-1/2}_{X_\text{can}}$ \cite{cosserat2023numerical}, we have $\Psi^*\Omega_\text{can} = \omega_\text{can}$, and the corresponding bi-realizations are
		\[\alpha_\text{can}=(Q+\frac{1}{2}\xi_P,P-\frac{1}{2}\xi_Q,Z+\frac{1}{2}\xi_T,T-\frac{1}{2}\xi_Z), \quad \beta_\text{can}(x,\xi) = \alpha_\text{can}(x,-\xi).\]
		These give an explicit, simple formula for a symplectic bi-realization in the canonical coordinates. Finally, transposing this symplectic bi-realization to the original homogeneous symplectic setting as $\alpha = F^{-1}\circ \alpha_\text{can} \circ (F^{-1})^*$ which result in 
		\begin{equation}\label{eq:contact bi-realization}
			\alpha(q,p,z,t,\xi_{q},\xi_{p},\xi_z,\xi_t)= \left(q-\frac{\xi_p}{2t},\ \frac{tp+\frac{\xi_q}{2}}{t-\frac{\xi_z}{2}}, z+\frac{\xi_t}{2} - \frac{p\xi_p}{2t}, t-\frac{\xi_z}{2}\right)
		\end{equation}
		with $\beta(x,\xi) = \alpha(x,-\xi)$.}
	
	Even in this contact canonical example, finding an exact bi-realization is already nontrivial, and for more complex Poisson or Jacobi structures it becomes nearly impossible. A practical workaround is to consider an approximate bi-realization, which is sufficient for our purpose of constructing numerical integrators. Following \cite{cabrera2024}, one can compute an approximation of Karasev's realization up to order 3 of Karasev's realization given by 
	\begin{align}\label{eq: approximation alpha}
		\alpha^i(y,\xi) = y^i - \frac{1}{2}\pi^{vi}\xi_v - \frac{1}{12}\partial_u\pi^{vi}\pi^{wu}\xi_v\xi_w -\frac{1}{48}\partial_u\partial_w \pi^{vi}\pi^{ku}\pi^{lw}\xi_v\xi_k\xi_w.
	\end{align}
	
	If we look closely, this approximation preserves the original homogeneity of $\alpha$. We have two cases:
	\begin{itemize}
		\item[(i)] $i\neq t$, component $0$-homogeneous: when all the variables $\xi$ are of the form $\xi_j$ with $j\neq t$, the component $\pi^{ji}$ is $\frac{\Lambda^{ji}}{t}$, so we have a $0$-homogeneous term, and when we have $\xi_t$, the components $\pi^{ji}=-E_j$, which are not homogeneous.
		
		\item[(ii)]$i=t,$ component $1$-homogeneous: is similar to the previous case, except that the $\pi^{vi}$ is always $-E_j$, that are multiplied by $\xi_j$ that are homogeneous. So in this case exits always a variable $1$-homogeneous.
	\end{itemize}
	So $\alpha^i\circ T^*h_z = h_z\circ \alpha^i.$
	
	In our case, the approximation of order $1$ of the realization is
	\begin{align}\label{eq: bi-realization example}
		\alpha(q,p,z,t,\xi_{q},\xi_{p},\xi_z,\xi_t) &= \left(q - \frac{1}{2t}\xi_{p}, p_i+\frac{1}{2t}\left(\xi_{q} + p\xi_z\right), z - \frac{1}{2}\left(\frac{p}{t}\xi_{p_i} - \xi_t\right), t-\frac{1}{2}\xi_z\right),
	\end{align}
	\begin{align}\label{eq: bi-realization example2}
		\beta(q,p,z,t,\xi_{q},\xi_{p},\xi_z,\xi_t) &=  \left(q + \frac{1}{2t}\xi_{p}, p-\frac{1}{2t}\left(\xi_{q} + p\xi_z\right), z + \frac{1}{2}\left(\frac{p}{t}\xi_{p} - \xi_t\right), t+\frac{1}{2}\xi_z\right).
	\end{align}
	
	{
		\begin{remark}
			As mentioned before, finding explicit homogeneous bi-realizations of homogeneous Poisson manifolds can be hard already when they are in canonical form. For other cases, it is to be expected that besides making use of approximate bi-realizations (or rather, before doing that), it will be helpful to use Weinstein's splitting theorem \cite{Weinstein1983} in order to obtain better adapted coordinates to the Poisson manifold at hand; when possible, it would be even better to put the Poisson manifold in canonical form, for example around a symplectic leaf \cite{CrainicMarcut2012} or a Poisson manifold \cite{FernandesMarcut2026}.
		\end{remark}
	}
	
	\section{Homogeneous symplectic groupoids for homogeneous integrators}\label{sec-smoothfamilies}
	
	This section is central not only for constructing the Jacobi-Hamiltonian integrator, but also for truly understanding how it works in practice. A rough outline of the construction, for a Hamiltonian vector field $X_H$ on a Jacobi manifold $J$, is as follows:
	\begin{enumerate}
		\item Poissonization: Lift $J$ to a Poisson manifold $\tau: P_J \to J$ and define a Hamiltonian vector field $X^P_H$ on $P_J$ inducing the same dynamics on $J$ as $X_H$.
		\item Homogeneous bi-realization: Construct an explicit homogeneous bi-realization of $P_J$. 
		\item Homogeneous Lagrangian bisections: Construct a family of particular homogeneous Lagrangian submanifolds (bisections) of the bi-realization, generating the flow of $X^P_H$.
		\item Flow approximation: Approximate the flow by approximating the family of homogeneous Lagrangian bisections. 
	\end{enumerate}
	\newpage
	
	{\begin{remark}[On constructivity]
			Step $1$ is fully constructive. Step $2$ is theoretically constructive, though exact implementation can be challenging; approximations can be used, as illustrated in Section~\ref{sec: example}. Step $3$ is studied in detail in this section. Step $4$ is not addressed here and will be explored in a subsequent work. 
	\end{remark}}

	In this section, we focus on step 3 of the construction outlined above. We introduce smooth families of Lagrangian submanifolds, examine their normal variations and the associated variation forms, and use these forms to highlight essential structural properties of the family.
	
	
	\subsection{Smooth families of homogeneous Lagrangian bisections}
	
	\begin{definition}[Smooth family]\label{def: L_k} 
		Let $(\Sigma,\omega_\Sigma,h)$ be a homogeneous symplectic manifold. A family $(L_s)_{s\in I}$ of $\mathbb{R}^\times$-submanifolds of $\Sigma$ parameterized by $I$ is said to be a smooth family of homogeneous Lagrangian submanifolds if all $L_s$ are homogeneous Lagrangian and $L_I =\{(x,s)\in \Sigma\times I, x\in L_s\}$ is a submanifold of $\Sigma\times I$ such that the restriction to $L_I$ of projection $S\times I\rightarrow I$ is a surjective submersion.
	\end{definition}
	
	From now on, fix a smooth family of homogeneous Lagrangian submanifolds $(L_s)_{s\in I}$ of $\Sigma$ as in Definition \ref{def: L_k}. Let $NL_s = T\Sigma|_{L_s}/TL_s$ be the normal bundle of $L_s$. Let us describe the construction from \cite{COSSERAT2023} of the section $\left[\frac{\partial L_s}{\partial s}\right]\in \Gamma(NL_{s_0})$, called the \textit{normal variation} of $(L_s)_{s\in I}$ at $s_0$. At a point $x\in L_{s_0}$, it is defined by $\left[\frac{\partial L_s}{\partial s}\right]:=\left[\frac{\partial \gamma(s)}{\partial s}_{|s=s_0}\right]\in (NL_s)_x$, where $\gamma: I\to \Sigma$ is any \textit{$L$-path through $x$}, meaning a smooth path such that $\gamma(s)\in L_s$ and $\gamma(s_0)=x$. Lemma 2.3 from \cite{COSSERAT2023} guarantees that the normal variation is well defined and smooth. In particular, its value at $x$ is independent of the choice of $L$-path through $x$. Since we are in a homogeneous symplectic manifold, $NL_s$ is canonically isomorphic to $T^*L_s$ by Theorem \ref{thm-WLtub} and the normal variation corresponds to a family of 1-homogeneous $1$-forms $\xi_s\in \Omega^1(L_s)$, called \textit{variation forms}, and satisfy the equation
	\begin{equation}\label{eq: variation form}
		\omega_\Sigma\left(\left[\frac{\partial L_s}{\partial s}(x)\right], u\right) =  \xi_s(u),\ \forall u\in T_xL_s.
	\end{equation}
	
	\begin{proposition} The normal variations $\left[\frac{\partial L_s}{\partial s}(x)\right]$ of a smooth family of homogeneous Lagrangian manifolds are 0-homogeneous. Equivalently, the variation forms $\xi_s$ are 1-homogeneous.
	\end{proposition}\label{prop: normal variations and variation forms}
	
	\noindent \textit{Proof:}
		Let $x\in L_{s_0}$, and let $\gamma$ be an $L$-path through $x$. Then, since each $L_s$ is homogeneous and $\gamma(s)\in L_s$, we know that $h_z(\gamma(s))\in L_s$, so $h_z\circ \gamma$ is an $L$-path through $h_z(x)$. Thus  
		
		\[h_{z*}\left[\frac{\partial L_s}{\partial s}(x)\right]=h_{z*}\left[\frac{\partial \gamma(s)}{\partial s}_{|s=s_0}\right]=\left[\frac{\partial (h_z\circ\gamma)(s)}{\partial s}_{|s=s_0}\right]=\left[\frac{\partial L_s}{\partial s}(h_z(x))\right],\]
		so the normal variation is 0-homogeneous.
		From Equation (\ref{eq: variation form}), we have that
		\begin{align*}
			(h_{z}^*\xi_s)(u) &=\left(h_{z}^*\left(i_{\frac{\partial L_s}{\partial s}}\omega_\Sigma\right)\right)(u)  \overset{(\ref{eq: relation})}{=} (h^*_{z}\omega_\Sigma)\left(h_{z*}^{-1}\left[\frac{\partial L_s}{\partial s}\right], u\right)=z\omega_\Sigma\left(h_{z*}^{-1}\left[\frac{\partial L_s}{\partial s}\right], u\right).
		\end{align*} 
		Since $z\xi_s(u)=z\omega_\Sigma\left(\left[\frac{\partial L_s}{\partial s}\right], u\right)$, we conclude that 1-homogeneity of the variation form $\xi_s$ is equivalent to the 0-homogeneity of $\left[\frac{\partial L_s}{\partial s}(x)\right]$. 
	
	\noindent $\square$

	\begin{definition}[Exact smooth family of Lagrangian submanifolds]
		We call exact a smooth family of Lagrangian submanifolds $(L_s)_{s\in I}$ such that its corresponding variation $1$-forms $(\xi_s)_{s\in I}$ are exact; in that case we call variation functions to any time-dependent functions $(f_s)_{s\in I}$ such that $df_s = \xi_s$, for all $s\in I$.
	\end{definition}
	
	\begin{corollary}
		Let $(L_s)_{s\in I}$ be an exact smooth family of Lagrangian submanifolds with variation forms $(\xi_s)_{s\in I}$ and variation functions $(f_s)_{s\in I}$. Then $\xi_s$ are $1$-homogeneous if and only if $f_s$ are $1$-homogeneous.
	\end{corollary}

	Now, we exhibit two examples of smooth families of homogeneous Lagrangian submanifolds which will be of use.
	\begin{example}\label{ex: dH}
		Let $H\in C^\infty(\Sigma)$ be a $1$-homogeneous Hamiltonian function whose $X_H\in\mathfrak{X}(\Sigma)$ is complete. Let $L$ be a homogeneous Lagrangian submanifold. The family $L_s = \phi^s_H(L)$ is an exact smooth family of homogeneous Lagrangian submanifolds, and the variation form at $t$ is $dH$.
	\end{example}
	\begin{example}\label{ex: cotangent projec}
		Let $T^*Q$ be the cotangent bundle of a $\mathbb{R}^\times$-manifold $Q$. For every family of closed homogeneous $1$-forms $(\zeta_s)_s$, their images $L_s = \{\zeta_s(x),\ x\in Q\}$ are a smooth family of homogeneous Lagrangian submanifolds, using Proposition \ref{prop: mu closed}. The homogeneous variation form at $t$ is $\tau^*\partial_s \zeta_s$, where $\tau$ is the cotangent projection.
	\end{example}
	
	\begin{proposition}\label{prop: two symplectic lagrangian}
		Let $(V,\omega_V,h^V)$ and $(W,\omega_W,h^W)$ be two homogeneous symplectic manifolds, $\phi:V\overset{\sim}{\rightarrow} W$ a homogeneous symplectomorphism, that is, $\phi\circ h^V = h^W\circ \phi$ and $(L_s)_s$ a smooth family of homogeneous Lagrangian submanifolds on $V$ with homogeneous variation forms $\xi_s$. Then, $\tilde{L}_s = \phi(L_s)$ is also a smooth family of homogeneous Lagrangian submanifolds with homogeneous variation forms $\tilde{\xi}_s$ such that $\xi_s = \phi^*\tilde{\xi_s}.$
	\end{proposition}
	
	\noindent \textit{Proof:}
		Since $\phi$ is a homogeneous symplectomorphism, we only need to prove that $\tilde{L_s}$ are Lagrangian submanifolds and that their variation forms are $\tilde{\xi_s}.$ Let $x\in T\tilde{L}_s$, then there exists $y\in TL_s$ such that $x=\phi_*y.$ So for any $\tilde{u}\in TW$,
		\begin{equation}
			\omega_W(x,\tilde{u}) = \omega_W(\phi_*y,\phi_*u) = (\phi^*\omega_W)(y,u) = \omega_V(y,u) =0.
		\end{equation}
		
		Now, suppose that $\tilde{\xi}_s$ are the homogeneous variation forms of $\tilde{L}_s$, so that they satisfy the relation
		\begin{align*}
			&\tilde{\xi}_s(\tilde{u}) = \omega_W\left(\left[\frac{\partial \tilde{L}_s}{\partial s}\right], \tilde{u}\right) = \ \omega_W\left(\phi_*\left[\frac{\partial L_s}{\partial s}\right], \tilde{u}\right) = \left(i_{\phi_*\left[\frac{\partial L_s}{\partial s}\right]} \omega_W \right)(\tilde{u})\\
			\Rightarrow & \phi^*\tilde{\xi}_s = \phi^*\left(i_{\phi_*\left[\frac{\partial L_s}{\partial s}\right]} \omega_W\right) \overset{(\ref{eq: relation})}{=} i_{\left[\frac{\partial L_s}{\partial s}\right]}\phi^*\omega_W = i_{\left[\frac{\partial L_s}{\partial s}\right]}\omega_V = \xi_s. 
		\end{align*}
	
	\noindent $\square$  
	
	
	\subsection{Homogeneous Hamilton-Jacobi equation}
	
	Let $Q$ be a $\mathbb{R}^\times$-manifold and $T^*Q$ its homogeneous cotangent bundle. Consider the $1$-homogeneous Hamiltonian function $H\in\Omega^0(Q)$ and the homogeneous Hamilton-Jacobi equation (\ref{eq: H-J}). There are two families of Lagrangian submanifolds related:
	\begin{example}
		The Hamiltonian flow $\phi^t_H: T^*Q\rightarrow T^*Q$ is a homogeneous symplectomorphism, $\phi^t_H\circ T^*h_z =  T^*h_z\circ\phi^t_H$. Using Proposition \ref{prop: phi lagrangian}, the graphs $G^t_H =\{(x,\phi^t_H(x))\in T^*Q\times T^*Q\}$ are homogeneous Lagrangian submanifolds of $T^*Q\times T^*Q$. The map $\Phi^{-t}_H: G^t_H\rightarrow T^*Q,\ (x,\phi^t_H(x))\mapsto \phi^t_H(x) $ is a homogeneous symplectomorphism. By Example \ref{ex: dH}, the variation form in $T^*Q$ is $dH$ so, using Proposition \ref{prop: two symplectic lagrangian}, the variation form of $G^t_H$ is $\Phi^{-t*}_HdH.$
	\end{example}
	
	\begin{example}
		Let $\mathbf{S}_t$ be the solution of Hamilton-Jacobi equation (\ref{eq: H-J}). It is a $1$-homogeneous function on $Q\times Q$ and its differentials $(d\mathbf{S}_t)_t$ are also a $1$-homogeneous exact and closed forms. Using Proposition \ref{prop: mu closed}, their images $\underline{d\mathbf{S}_t}$ are exact homogeneous Lagrangian submanifolds of $T^*(Q\times Q)$. By Example \ref{ex: cotangent projec}, their variation forms are $\tau^*d\frac{\partial \mathbf{S}_t}{\partial t}$.
	\end{example}
	
	These two variation forms are related by the homogeneous symplectomorphism
	\begin{align}\label{eq: Psi}
		\begin{split}
			\Psi:T^*Q\times T^*Q&\rightarrow T^*(Q\times Q)\\
			(\xi(q),\bar{\xi}(\bar{q}))&\mapsto \xi(q)-\bar{\xi}(\bar{q}).
		\end{split}
	\end{align}
	If we define $T^*h^\times_z = (T^*h\times T^*h)_z$ the action in $T^*Q\times T^*Q$ and $T^*h^Q_z$ the lifted action of $h^Q_z = (h\times h)_z$ in $Q\times Q$, we can see homogeneity
	\begin{align*}
		(\Psi\circ T^*h^\times_z)(\xi(q),\bar{\xi}(\bar{q})) &=\Psi(T^*h_z\xi(q),T^*h_z\bar{\xi}(\bar{q})) = T^*h_z\xi(q) - T^*h_z\bar{\xi}(\bar{q}) \\
		&= T^*h^Q_z(\xi(q)-\bar{\xi}(\bar{q})) = (T^*h^Q_z\circ \Psi)(\xi(q),\bar{\xi}(\bar{q})).
	\end{align*}
	
	
	\subsection{Homogeneous symplectic groupoids}
	
	We have seen with Proposition \ref{thm: bi-realization} that for a homogeneous Poisson manifold, we can construct a homogeneous symplectic bi-realization.
	To continue the construction, let us define groupoids and symplectic groupoids (see \cite{crainic2021lectures} for a textbook account).
	\begin{definition}[Groupoid and Lie groupoid]
		A groupoid, denoted as $\mathcal{G}\rightrightarrows M$, is a set $M$ of objects and a set $\mathcal{G}$ of arrows, together with the following structure maps:
		\begin{enumerate}
			\item[(i)] source $s:\mathcal{G}\rightarrow M$ and target  $t:\mathcal{G}\rightarrow M$;
			\item[(ii)] multiplication $\textbf{m}:\mathcal{G}^{(2)}\rightarrow \mathcal{G},\ (g,h)\mapsto \textbf{m}(g,h) :=g\cdot h$, where
			\begin{align*}
				\mathcal{G}^{(2)} := \{(g,h)\in \mathcal{G}\times\mathcal{G}: s(g)=t(h)\}
			\end{align*}
			and which satisfies
			\begin{enumerate}
				\item[-]  $s(g\cdot h) = s(h)$ and $t(g\cdot h) = t(g)$,
				\item[-] $(g\cdot h)\cdot k = g \cdot( h\cdot k)$;
			\end{enumerate}
			\item[(iii)] unit map $\sigma:M\rightarrow \mathcal{G},\ x\mapsto \sigma(x) := 1_x$ which satisfies
			\begin{enumerate}
				\item[-] $s(1_x) = t(1_x) = x$,
				\item[-] $g\cdot 1_{s(g)} = 1_{t(g)}\cdot g =g$;
			\end{enumerate}
			\item[(iv)] inverse map $i:\mathcal{G}\rightarrow \mathcal{G},\ g\mapsto i(g) :=g^{-1}$, which satisfies
			\begin{enumerate}
				\item[-] $s(g^{-1}) = t(g)$ and  $t(g^{-1}) = s(g)$,
				\item[-] $g^{-1}\cdot g = 1_{s(g)}$ and $g\cdot g^{-1} = 1_{t(g)}$. 
			\end{enumerate}
		\end{enumerate}
		
		If $\mathcal{G},\ M$ are manifolds, $s,t$ are submersions, and $\textbf{m},\sigma$ and $i$ are smooth maps, we say that $\mathcal{G}\rightrightarrows M$ is a Lie groupoid.	
	\end{definition}
	
	\begin{definition}[Symplectic groupoid]
		A symplectic groupoid is a Lie groupoid $\Sigma\rightrightarrows M$ with a symplectic form $\omega\in \Omega^2(\Sigma)$ that satisfies the property
		\begin{align*}
			\textbf{m}^*\omega = pr_1^*\omega + pr_2^*\omega,
		\end{align*}
		where $pr_1,\ pr_2:\mathcal{G}^{(2)}\rightarrow \mathcal{G}$ are the projections of the first and second components.  
	\end{definition}
	
	\begin{definition}[Bisection]
		Let $\mathcal{G}\rightrightarrows M$ be a Lie groupoid. A submanifold $L\subset \mathcal{G}$ is called a bisection if the restrictions of both source and target to $L$ are diffeomorphisms onto $M$.
	\end{definition}
	
	\begin{remark}
		Every symplectic groupoid induces a symplectic bi-realization. In particular, it induces a Poisson structure on its base.    
		Not every Poisson manifold $(P,\pi)$ admits a symplectic groupoid which induces the Poisson structure $\pi$ on $P$. Nonetheless, any $(P,\pi)$ is induced by a \emph{local symplectic groupoid}, which is still a bi-realization; we say that the local symplectic groupoid integrates the Poisson manifold. These can be constructed via a spray, starting with the construction of bi-realization that we have described, as done in \cite{CMS2022}. 
	\end{remark}
	
	Denote by $\alpha$ the homogeneous source and by $\beta$ the homogeneous target of the previously constructed bi-realization.
	
	\begin{remark}
		We know that $\alpha$ is homogeneous, that is, $h_z\circ \alpha = \alpha\circ T^*h_z$. Moreover, given a homogeneous bisection $L$, the inverse restricted to $L$, $\alpha^{-1}_{|L}$ is also homogeneous: $\alpha^{-1}_{|L}\circ h_z = T^*h_z\circ \alpha^{-1}_{|L}$. And so, any bisection induces a homogeneous diffeomorphism of the unit manifold $M$ by $\phi_L\coloneqq \beta\circ \alpha^{-1}_{|L},$ because
		\begin{align*}
			\phi_L\circ h_z = (\beta\circ \alpha^{-1}_{|L})\circ h_z = \beta\circ (T^*h_z\circ\alpha^{-1}_{|L}) = h_z\circ (\beta\circ \alpha^{-1}_{|L}) = h_z\circ \phi_L.
		\end{align*}
	\end{remark}
	
	So, we have the following proposition.
	\begin{proposition}[\cite{CosteDazordWeinstein}]
		Let $(P,\Pi)$ be a homogeneous Poisson manifold with action $h$ and let $(\Sigma\rightrightarrows P_J,\Omega, T^*h)$ be a local symplectic groupoid integrating it. If a bisection $L\subset \Sigma$ is { homogeneous} Lagrangian, then:
		\begin{enumerate}
			\item[(i)] the induced diffeomorphism $\phi_L:P_J\rightarrow P_J$ is a homogeneous Poisson diffeomorphism;
			\item[(ii)] provided that the fibers of the source map are connected for all $x\in P$, $\phi_L(x)$ and $x$ belongs to the same symplectic leaf of $P$.
		\end{enumerate}
	\end{proposition}
	\newpage
	
	The Jacobi integrators we construct will be related with this proposition. Now, we are interested in smooth families of homogeneous Lagrangian bisections of a homogeneous symplectic groupoid $\Sigma$, parametrized by $I\subset \mathbb{R}$ an interval containing 0.
	
	\begin{example}
		Let $(\phi_t)$ be a smooth family of homogeneous symplectomorphisms of $(\Sigma, \omega_\Sigma)$, a  homogeneous symplectic manifold with action $h_z$. This family is the flow of a time-dependent $0$-homogeneous vector field $\overset{\rightarrow}{\xi}_t$ related by $\omega_\Sigma$ with a time-dependent $1$-homogeneous closed form $(\zeta)_t$, that is, $\overset{\rightarrow}{\xi}_t = \omega_\Sigma^{-1}(\zeta_t)$.
		Consider the pair groupoid $\Sigma\times\Sigma\rightrightarrows \Sigma$, equipped with the symplectic form $\Omega=pr_1^*\omega_\Sigma - pr_2^*\omega_\Sigma$. It is a symplectic groupoid over $(\Sigma, \omega_\Sigma)$.
		
		Any smooth family of homogeneous Lagrangian bisection $(L_\epsilon)_{\epsilon\in I}$ of $\Sigma\times\Sigma$ will be based on the choice of the first and second factors in $\Sigma\times\Sigma$ of the form $\{(x,\phi_\epsilon(x),\ x\in \Sigma\}_{\epsilon\in I}$.
		
		For instance, for any solution $\mathbf{S}_t$ of (\ref{eq: H-J}), a smooth family of Lagrangian bisections of the pair groupoid is given by $\Psi^{-1}(\underline{d\mathbf{S}_t}) = \{(d_q\mathbf{S}_t(q,\bar{q}), -d_{\bar{q}}\mathbf{S}_t(q,\bar{q})),\ (q,\bar{q})\in Q\times Q\}\subset T^*Q\times T^*Q$, where $\Psi$ is given by (\ref{eq: Psi}).
	\end{example}
	
	\begin{remark}
		Any exact family of homogeneous Lagrangian bisections $(L_t)_t$ naturally induces a homogeneous Poisson Hamilton Integrator (to be defined in the next Section) with time step $\Delta t$ by:
		\begin{align*}
			P_J&\rightarrow P_J\\
			x &\mapsto \beta\circ (\alpha_{|L_{\Delta t}})^{-1}(x). 
		\end{align*}
	\end{remark}
	
	
	\section{Jacobi Hamiltonian integrators}\label{sec: JHI}
	
	In this section, we define Jacobi Hamiltonian Integrators (JHI). Let us start with a Hamiltonian system on a Jacobi manifold $(J, \Lambda, E)$ defined by
	\begin{align*}
		X_H = \Lambda(\cdot,dH) - HE(\cdot),
	\end{align*}
	where $H$ is the Hamiltonian function. 
	
	\begin{definition}[Jacobi Hamiltonian integrator]
		Let $H^J\in C^\infty(J)$ be a Hamiltonian on $J$. A smooth family of diffeomorphisms of $J$, $(\phi_\epsilon)_\epsilon$ is a  Jacobi Hamiltonian integrator of order $k\geq 1$ for $H^J$ if:
		\begin{enumerate}
			\item $\phi_\epsilon$ is a Jacobi diffeomorphism\footnote{Let $(J_1,\Lambda_1,E_1)$ and $(J_2,\Lambda_2,E_2)$ be two Jacobi manifolds, $\phi:J_1\rightarrow J_2$ a diffeomorphism and $\{\cdot,\cdot\}_{J_i},\ i=1,2$ the respective Jacobi brackets. We say that $\phi$ is a Jacobi diffeomorphism if $\phi\circ \{f,g\}_{J_1} = \{\phi\circ f,\phi\circ g\}_{J_2}.$};
			\item there exists $(\mathcal{H}_s)_s$ a time-dependent Hamiltonian such that
			\begin{enumerate}
				\item $\mathcal{H}_s = H^J + o(s^{k-1})$,
				\item $\phi_\epsilon = \Phi^\epsilon_{(\mathcal{H}_s)_s}$ its the time-$\epsilon$ flow of $\mathcal{H}_s$.
			\end{enumerate}
		\end{enumerate}   
	\end{definition}
	
	To construct such an integrator, we first transform the Jacobi manifold into a homogeneous Poisson manifold $(P,\Pi,Z)$ as in Section \ref{sec: poissonization}. Then we can transport the original Hamilton function $H(x)$ to the homogeneous Poisson manifold considering $H^P(x,t) = tH(x)$. This will be a $1$-homogeneous Hamiltonian, so we can define a homogeneous Poisson-Hamiltonian system $X^P_H = \Pi(dH^P)$.
	
	Now, we want to construct a homogeneous Poisson Hamilton Integrator (hPHI) for $X^P_H$.
	
	\begin{definition}[Homogeneous Poisson Hamiltonian integrator] Let $(P,\Pi,h)$ be a homogeneous Poisson manifold and let $H^P\in C^\infty(P)$ be a $1$-homogeneous Hamiltonian on $P$. A smooth family of homogeneous diffeomorphisms of $P$, $(\phi_\epsilon)_\epsilon$ is a homogeneous Poisson Hamilton integrator of order $k\geq 1$ for $H^P$ if:
		\begin{enumerate}
			\item $\phi_\epsilon$ is a homogeneous Poisson diffeomorphism;
			\item there exists $(\mathcal{H}_s)_s$ a time-dependent homogeneous Hamiltonian such that
			\begin{enumerate}
				\item $\mathcal{H}_s = H^P + o(s^{k-1})$,
				\item $\phi_\epsilon = \Phi^\epsilon_{(\mathcal{H}_s)_s}$ is the time-$\epsilon$ flow of $\mathcal{H}_s$.
			\end{enumerate}
		\end{enumerate}   
	\end{definition}
	\newpage
	
	Using the following theorem, we can construct a homogeneous Poisson Hamilton Integrator by leveraging the families of exact Lagrangian bisections constructed previously.
	
	\begin{thm}[\cite{cabrera2024}]\label{thm: ham flow}
		Let $R=(\Sigma,\omega,\alpha,\beta,\sigma)$ be a symplectic bi-realization for $(P,\Pi)$ a Poisson manifold, where $\sigma$ is the unit map of a local symplectic structure on the bi-realization. Then;
		\begin{enumerate}
			\item when $L \xhookrightarrow{} (\Sigma,\omega)$ is a Lagrangian bisection for $R$, the induced map $\varphi_L = \beta\circ \alpha^{-1}_{|L}$ defines a Poisson diffeomorphism $(P,\Pi)\rightarrow (P,\Pi)$;
			\item if $\phi^s_H:P\rightarrow P$ is the Hamiltonian flow on $(P,\Pi)$ defined by the Hamiltonian function $H$, then
			\[\phi^s_H = \alpha\circ \phi^s_{\alpha^*H}\circ \sigma,\]
			with $\phi^t_{\alpha^*H}:\Sigma\rightarrow\Sigma$ the Hamiltonian flow of $\alpha^*H$ in $(\Sigma,\omega)$;
			\item in the previous item,
			\[\phi^s_H = \varphi_{L_s},\ \textnormal{for the Lagrangian bisection } L_s = \phi^s_{\alpha^*H}(\sigma(P)).\]
		\end{enumerate}
	\end{thm}
	
	In our case, this theorem can be applied because the unit map $\sigma$ for the symplectic bi-realizations that we constructed is the zero section of $T^*P$ \cite{Weinstein1971}, so it is homogeneous:
	\[T^*h_z\circ \sigma = \sigma \circ h_z.\]
	With this, as $H^P$ is $1$-homogeneous, $\alpha^*H^P$ is also $1$-homogeneous \[(T^*h_z)^*\alpha^*H^P = \alpha^*h_z^*H^P = \alpha^*(zH^P) = z\alpha^*H^P.\] So its flow is homogeneous and both sides of the item $2.$ agree on homogeneity.
	
	Theorem \ref{thm: ham flow} guarantees that, given a Lagrangian bisection $L$, the induced diffeomorphism $\varphi_L$ is a Hamiltonian flow.
	Given an exact smooth family of Lagrangian bisections such that $L_0$ is the zero section, constructing a suitable approximation of the variation functions $(\mathcal{H}_s)_s$ that coincides with $ H^P$ of order $k$, the induced family of diffeomorphisms $(\varphi_{L_s})_{s\in I}$ is a homogeneous Hamiltonian Poisson integrator of order $k$ for $H^P$.
	
	Using the equivalence between the Jacobi category and homogeneous Poisson category (Proposition \ref{Prop-equiv-cat}), we can conclude that we have a $1$-to-$1$ correspondence between homogeneous Poisson Hamiltonian integrators and Jacobi Hamiltonian integrators. This grants the existence of JHI's, constructed as explained above.
	
	In summary, the combined use of the Poissonization procedure, homogeneous symplectic bi-realizations, and smooth families of homogeneous Lagrangian bisections,  provides a constructive and geometrically natural method to produce structure-preserving numerical integrators for Jacobi Hamiltonian systems, extending the theory of Poisson Hamiltonian integrators to the Jacobi setting.
	
	
	{
		\subsection{Numerical example}
		To illustrate the Jacobi Hamiltonian integrator, we consider the damped harmonic oscillator, which can be formulated as a dissipative contact Hamiltonian system. We take the canonical contact structure introduced in Section~\ref{sec: example}, given by
		\[\Lambda = \frac{\partial}{\partial p}\wedge\frac{\partial}{\partial q} + p\frac{\partial}{\partial p_i}\wedge\frac{\partial}{\partial z},\qquad E=-\frac{\partial}{\partial z},\] and the Hamiltonian \[H(q,p,z) = \frac{p^2}{2} + \frac{q^2}{2} + \gamma z,\] where $\gamma\in\mathbb{R}$ is de damped parameter.
		Through Poissonization (symplectization), we obtain the Poisson structure \[\Pi = \frac{1}{t}\frac{\partial}{\partial p}\wedge\frac{\partial}{\partial q} + \left(\frac{p}{t} \frac{\partial}{\partial p} - \frac{\partial}{\partial t} \right)\wedge \frac{\partial}{\partial z} \] and the associated homogeneous Hamiltonian
		\[\hat{H}(q,p,z,t) = tH(q,p,t),\]
		which induces the original dynamics. Using the homogeneous symplectic bi-realization computed in \eqref{eq:contact bi-realization}, a first-order Jacobi Hamiltonian integrator is defined by the scheme:
		\begin{enumerate}
			\item solve for $y_n$,
			\[
			\alpha\big(y_n, \Delta s \nabla \hat{H}(y_n)\big) = x_n;
			\]
			\item update the solution,
			\[
			x_{n+1} = \beta\big(y_n, \Delta s \nabla \hat{H}(y_n)\big).
			\]
		\end{enumerate}
		On $[0,20]$, with step size $\Delta s = 0.5$, initial condition $x_0 = (1,0,0)$, and $\gamma = 0.05$, Figure~\ref{fig:contactdamped} shows that the Jacobi Hamiltonian integrator more accurately captures the dissipative dynamics than the semi-implicit symplectic Euler method. The scheme is implemented directly using the explicit bi-realization \eqref{eq:contact bi-realization}, which allows the integrator to be evaluated efficiently at each step. 
		\begin{figure}[ht]
			\includegraphics[width=\textwidth]{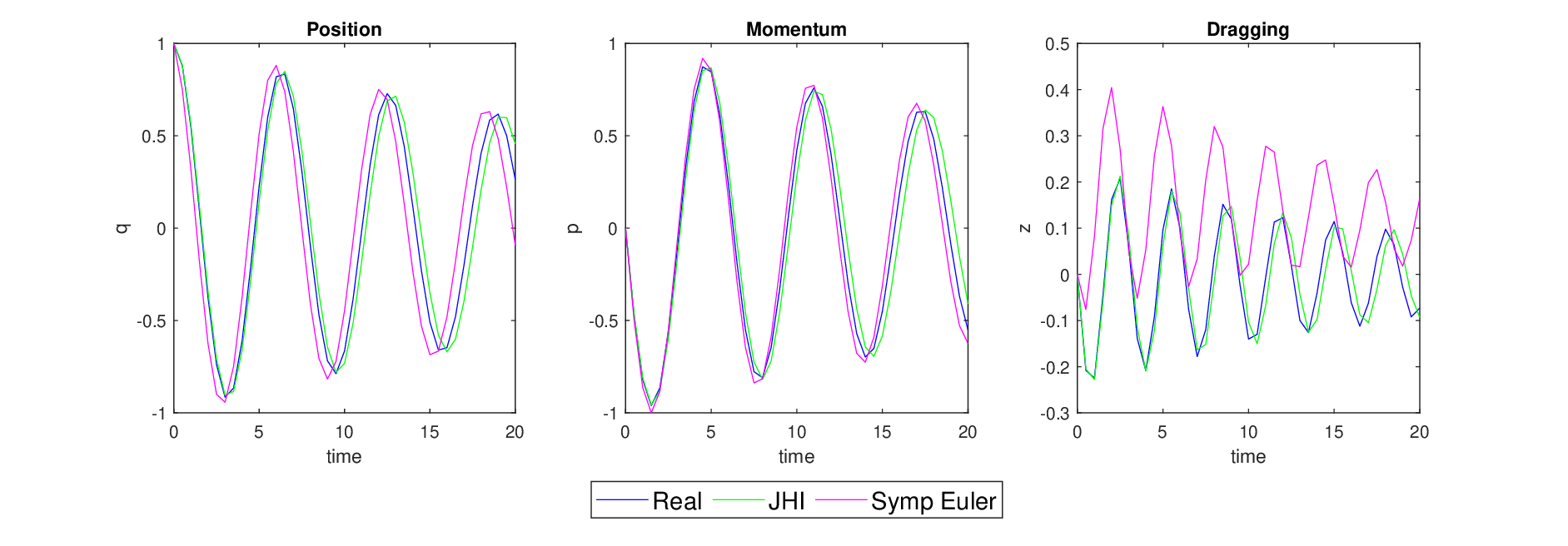}
			\caption{Trajectories of the damped parametric oscillator using first order JHI and symplectic Euler method.}
			\label{fig:contactdamped}
	\end{figure}}
	
	
	\section{Conclusion}
	
	This work has introduced Jacobi Hamiltonian Integrators, a new class of structure-preserving numerical schemes for Hamiltonian systems defined on Jacobi manifolds. The construction is based on lifting the problem to a homogeneous Poisson manifold via Poissonization, applying Poisson integrators in that setting, and projecting the result back to the original Jacobi manifold. This approach preserves not only the Jacobi structure and the Hamiltonian, but also the induced foliation and Casimir functions, which is particularly relevant in mechanical applications.
	
	A key ingredient in this construction is the interplay between contact, symplectic, and Poisson geometry. In particular, we make essential use of homogeneous symplectic realizations and homogeneous Lagrangian bisections. Under suitable conditions, we show that explicit symplectic realizations in the homogeneous Poisson context preserve the homogeneity, enabling the construction of homogeneous symplectic bi-realizations. These make use of homogeneous symplectomorphisms between a local symplectic groupoid $\Sigma \rightrightarrows P$ and a neighborhood of the zero section in $T^*P$, and allow for the transformation of solutions of the Hamilton-Jacobi equation into smooth families of homogeneous Lagrangian bisections $(L_t)_t$. Each such bisection determines, via composition with source and target maps, a homogeneous Poisson diffeomorphism $P\to P$ which corresponds to a Jacobi diffeomorphism. The resulting JHI method can thus be interpreted as a numerical approximation built from these structure-preserving transformations.
	
	Future work may focus on the development of explicit low-order JHI schemes and their analysis on concrete examples. Backward error analysis in the Jacobi setting could help clarify long-term behavior, while further exploration of contact groupoids and their discretizations might provide a path toward global integration methods on more general manifolds. On the computational side, efficient implementation and benchmarking against existing Poisson and symplectic schemes remain important challenges.
	
	Overall, this methodology extends the scope of geometric integration to Jacobi manifolds, contributing to a unified framework for the study of Hamiltonian and non-Hamiltonian systems with underlying geometric structure.
	
	
	\subsubsection*{Conflict of interest} The authors have no conflict of interest to disclose.
	
	\subsubsection*{Data Availability} There is no data associated to this work.
	
	\appendix
	
	\newpage%
	\renewcommand{\thesection}{\Alph{section}}
	
	
	\section{Normal forms in homogeneous symplectic geometry}\label{anexo: symplectic}
	
	{
		In this Appendix, we establish homogeneous versions of classical normal form results from symplectic geometry: the Poincar\'e Lemma~\ref{lemma-hPL}, the Darboux-Weinstein Theorem~\ref{APP-thm-dw}, the Weinstein Lagrangian neighborhood Theorem~\ref{thm: hWL}, and the Weinstein tubular neighborhood Theorem~\ref{App-thm-WLtub}. We briefly indicate which results follow directly from homogeneity and which require additional care.
		
		The homogeneous Poincar\'e Lemma~\ref{lemma-hPL} admits a simple, even more straightforward proof than in the classical case. Although more general versions exist (e.g., \cite[Lemma 10.3]{GG2025}), the form presented here suffices for our purposes.
		
		For the remaining local normal form results, it is crucial that all neighborhoods where they are defined are homogeneous (invariant under the $\mathbb{R}^\times$-action). This ensures that the constructions remain within the category of objects under study, so that the homogeneous diffeomorphisms correspond to isomorphisms of $\mathbb{R}^\times$-principal bundles.
		
		The proof of the homogeneous Darboux-Weinstein Theorem~\ref{APP-thm-dw} follows the standard argument, with neighborhoods defined via the flow of a homogeneous vector field. For the homogeneous Weinstein Lagrangian and tubular neighborhood Theorems~\ref{thm: hWL} and~\ref{App-thm-WLtub}, homogeneity of neighborhoods can be guaranteed using $\mathbb{R}^\times$-invariant Riemannian metrics and partitions of unity. Constructing the required homogeneous diffeomorphisms requires care: each auxiliary choice and identification in the standard proofs must respect the $\mathbb{R}^\times$-action. This is achieved as a consequence of $L$ being a homogeneous Lagrangian submanifold together with invariant metrics.
		
		These results provide the rigorous framework supporting the constructions in the main text, ensuring that all geometric objects and integrators are compatible with the $\mathbb{R}^\times$-structure intrinsic to Jacobi manifolds.}

	{ \begin{lemma}[Homogeneous Poincar\'e lemma]\label{lemma-hPL}
			Let $(M,h)$ be a $\mathbb{R}^\times$- manifold. Then all $1$-homogeneous closed $p$-forms are exact and admit $1$-homogeneous primitives.
		\end{lemma}
		\noindent \textit{Proof:}
			Let $\omega$ be a 1-homogeneous $p$-form and let $Z$ denote the infinitesimal generator of the principal $\mathbb{R}^\times$-action $h$. By definition of homogeneity,  
			\[
			\mathcal{L}_Z \omega = \omega.
			\]
			Using Cartan's magic formula and the closedness of $\omega$, we have  
			\[
			\mathcal{L}_Z \omega = i_Z d\omega + d(i_Z \omega) = d(i_Z \omega).
			\]
			Defining $\alpha = i_Z \omega$, it follows that  $
			\omega = d\alpha.$
			To check the homogeneity of $\alpha$, we compute $h_z^* \alpha $ for any $z \in \mathbb{R}^\times$, and using (\ref{eq: relation}) together with the 0-homogeneity of $Z$, we get
			\[
			h_z^* \alpha = h_z^* d(i_Z \omega) =  d(h_z^* (i_Z \omega))= d(i_Z h_z^* \omega) = d(i_Z z \omega) = z  \alpha.
			\]
			Hence, $\alpha$ is 1-homogeneous.
	
	\noindent $\square$ }
	
	With this, we can use a homogeneous version of Moser's trick to prove the following theorem.
	
	\begin{thm}[Homogeneous Darboux-Weinstein theorem]\label{APP-thm-dw}
		Let $(M,h)$ be a $\mathbb{R}^\times$-manifold and let $X\subset M$ be a $\mathbb{R}^\times$-submanifold. Suppose $\omega_0,\ \omega_1$ are two $1$-homogeneous symplectic forms on $M$, for which $(\omega_0)_{|X} =(\omega_1)_{|X} $. Then, there is a {homogeneous} neighborhood $\mathcal{U}$ of $X$ and a diffeomorphism $f:\mathcal{U}\rightarrow \mathcal{U}$ such that
		\begin{enumerate}
			\item $f(x) = x$, for all $x\in X$
			\item $f^*\omega_1 = \omega_0$
			\item if $h_z$ is the principal action, then $f\circ h = h\circ f.$
		\end{enumerate}
	\end{thm}
	\noindent \textit{Proof:}
		The proof is an adaptation of the proof of Theorem 3.2 of \cite{dwivedi2019hamiltonian}.
		
		Consider $\omega_s = (1-s)\omega_0 + s\omega_1$. For all $s\in[0,1]$, $\omega_s$ is closed since both $\omega_0$ and $\omega_1$ are closed. Since $d(\omega_0-\omega_1) =0$, we     can find a homogeneous $1$-form $\alpha$ such that $d\alpha =\omega_0-\omega_1 $, by the homogeneous Poincar\'e Lemma \ref{lemma-hPL}. Namely, $\alpha=i_Z(\omega_0-\omega_1)$, where $Z$ is the infinitesimal generator of the $\mathbb{R}^\times$-action.
		
		The hypothesis $(\omega_0)_{|X} =(\omega_1)_{|X}$, implies that $\alpha_{|X}=0$. Since $\omega_s|_X$ is symplectic for all $s\in[0,1]$, this is true for a small homogeneous neighborhood of $X$. Then we can find a time-dependent vector field $\eta_s$ such that
		\begin{align*}
			i_{\eta_s}\omega_s = \alpha.
		\end{align*}
		Note that, using $(\ref{eq: relation})$ we conclude that $\eta_s$ is $0$-homogeneous. So integrating $\eta_s$ gives us a family of local diffeomorphisms {(with $\mathbb{R}^\times$-invariant domains of definition, because they are the flow of a 0-homogeneous vector field)} $f_s$ with $f_0=id$, which commute with the action $h$, that is, $f_s\circ h = h\circ f_s$ and
		\begin{align*}
			\frac{d}{ds}f_s(m) = \eta_s(f_s(m)).
		\end{align*}
		We have also $(\eta_s)|_X=0$ and so $(f_s)|_X = id$. Using Proposition 6.4 in \cite{da2001lectures}, we have the following
		\begin{align*}
			(f_1)^*\omega_1 - \omega_0 & = \int_0^1 \frac{d}{ds}(f_s^*\omega_s)ds\\
			& = \int_0^1 f_s^*d(i_{\eta_s}\omega_s)ds + \int_0^1 f_s^*(\omega_1-\omega_0)ds\\
			& = \int_0^1 f_s^*d(\beta)ds + \int_0^1 f_s^*(\omega_1-\omega_0)ds\\
			& =\int_0^1 f_t^*(\omega_0-\omega_1)ds + \int_0^1 f_t^*(\omega_1-\omega_0)ds\\
			& =0.
		\end{align*}
		Thus, and because the domain of definition of $f_1$ is $\mathbb{R}^\times$-invariant, $f_1$ is the desired diffeomorphism. 
	
	\noindent $\square$

	We now present three auxilliary neighborhood results.
	
	{ 
		\begin{thm}[Homogeneous tubular neighborhood theorem]\label{thm: tubular nbh}
			Let $(M,h)$ be an  $n$-dimensional $\mathbb{R}^\times$-manifold and let $X$ be a $k$-dimensional $\mathbb{R}^\times$-submanifold, with $i:X\xhookrightarrow{} M$ the inclusion map, and $i_0:X\xhookrightarrow{} NX$ the embedding as the zero section of the normal bundle.
			
			Then, there exist a homogeneous and convex neighborhood $\mathcal{U}_0$ of $X$ in $NX$, a homogeneous neighborhood $\mathcal{U}$ of $X$ in $M$, and a homogeneous diffeomorphism $\varphi:\mathcal{U}_0\rightarrow \mathcal{U}$ such that the following diagram commutes
			\begin{center}
				\begin{tikzpicture}
					\node at (0,0) (a) {$\mathcal{U}_0$};
					\node at (1,-1) (b) {$X$};
					\node at (2,0) (c) {$\mathcal{U}$};
					\draw [->] (b) to node[left] {$i_0$} (a);
					\draw [->] (a) to node[above] {$\varphi$} (c);
					\draw [->] (b) to node[right] {$i$}(c);
				\end{tikzpicture}
			\end{center}
		\end{thm}
		
		\noindent \textit{Proof:} See the proof of Theorem \ref{thm: epsilon nbh}.
		
		\noindent $\square$ 
	}
	
	\begin{thm}[Homogeneous $\epsilon$-neighborhood theorem]\label{thm: epsilon nbh}
		Let $(M,h)$ be an  $n$-dimensional $\mathbb{R}^\times$-manifold and let $X$ be a $k$-dimensional $\mathbb{R}^\times$-submanifold. There exists a 0 - homogeneous function $\epsilon:X\to\mathbb{R}_+$ such that every point $p$ in the set $\mathcal{U}^\epsilon = \{p\in M\ |\ d(p,q)<\epsilon(q), \ \text{for some }q \in X\}$ has a unique nearest point $q\in X$ (where d denotes the Riemannian distance with respect to an $\mathbb{R}^\times$-invariant metric). 
		
		Moreover, setting $q=\pi(p)$, the map  $\pi:\mathcal{U}^\epsilon\to X$ is a submersion, and for every $p\in \mathcal{U}^\epsilon$ there is a unique geodesic curve joining $p$ to $q=\pi(p)$.
	\end{thm}
	{
		\noindent \textit{Proof:}
			The proofs of both these theorems follow closely those of their standard versions, Theorems 6.5 and Theorem 6.6 of \cite{da2001lectures} (see also \cite{lee-Riemannian}). Only two adaptations are needed to the homogeneous setting. The first is to choose a Riemannian metric $g$ on $M$ which is invariant under the $\mathbb{R}^\times$-action (that is, $\mathbb{R}^\times$ acts by isometries of $g$). It is always possible to find such a $g$ because the action is proper.
			
			That $g$ is $\mathbb{R}^\times$-invariant implies three properties: that also the Riemannian distance $d$ will be $\mathbb{R}^\times$-invariant; that the open neighborhoods of $X$ 
			
			\[\mathcal{U}^\epsilon = \{p\in M\ |\ d(p,q)<\epsilon(q), \ \text{for some }q \in X\}\subset M\] 
			and        
			\[NX^\epsilon = \{(x,v)\in NX\ | \ \sqrt{g_x(v,v)}<\epsilon(x)\}\subset T^* X\] are homogeneous, as long as $\epsilon$ is a $0$-homogeneous function; that the exponential map 
			
			\[exp:NX^\epsilon\rightarrow M,\] defined by $exp(x,v) = \gamma(1)$, where $\gamma:[0,1]\rightarrow M$ is the geodesic with $\gamma(0) = x$ and $\frac{d\gamma}{ds}(0) = v$, which maps $NX^\epsilon$ diffeomorphically to $\mathcal{U}^\epsilon$, is a homogeneous map.
			
			The other adaptation needed is to use $\mathbb{R}^\times$-invariant partitions of unity in order to guarantee the existence of a function $\epsilon$ as in the statement which is $0$-homogeneous, Once again, $\mathbb{R}^\times$-invariant partitions of unity exist because the action of $\mathbb{R}^\times$ is proper.
		
		\noindent $\square$ 
	}

	\begin{thm}[Homogeneous Whitney extension theorem]\label{thm: whitney extension}
		Let $(M,h)$ be an $n$-dimensional $\mathbb{R}^\times$-manifold and let $X$ be a $k$-dimensional $\mathbb{R}^\times$-submanifold. Suppose that at each $p\in X$ we have a linear isomorphism $L_p:T_pM\overset{\simeq}{\rightarrow}T_pM$ such that $L_{p|T_pX} = id_{T_pX}$, $L_p$ depends smoothly on $p$ and is homogeneous, that is $Th_z\circ L_p = L_{h_z(p)}\circ \dot{h_z}$ where $\dot{h_z}$ is the action of the fibers of $T_pM$. Then, there exists a homogeneous embedding $f:\mathcal{U}\rightarrow M$ of some homogeneous neighborhood $\mathcal{U}$ of $X$ in $M$ such that $f_{|X}=id_X$, $df_p = L_p$ for all $p\in X$.    
	\end{thm}
	
	\noindent \textit{Proof:}
		Let $g$ be a Riemannian metric on $M$ which is invariant ($0$-homogeneous) with respect to the action of $\mathbb{R}^\times$, and denote by $d$ its Riemannian distance (which is then also invariant under the action). Let \[\mathcal{U}^\epsilon = \{p\in M\ |\ d(p,X)\leq \epsilon\}\] be a neighbourhood of $X$ in $M$ for a $0$-homogeneous function $\epsilon:X\rightarrow \mathbb{R^+}$ which is small enough that any $p\in \mathcal{U^\epsilon}$ has a unique nearest point in $X$; define $\pi:\mathcal{U^\epsilon}\rightarrow X$, $p\mapsto$ nearest point to $p$ in $X$. If $\pi(p) = q$, then $p=exp(q,v)(1)$ for some $v\in N_qX = (T_qX)^\perp$.
		With assumptions, $\mathcal{U^\epsilon}$ is homogeneous: let $p\in \mathcal{U}^\epsilon$,
		\begin{align*}
			d(h_zp,X) = \underset{x\in X}{\text{inf}} d(h_zp,x) = \underset{x\in X}{\text{inf}} d(p,h_{z^{-1}}x) = d(p,X).
		\end{align*}
		
		Let $(x,t)$ be coordinates on $M$ with respect to a local trivialization, and let  $\dot{h}$ be the action of $Th$ only on the tangent fibers, for every $v\in T_qM$, $\dot{h_z}v = (v_x,zv_t)$. With this, we can prove that the exponential map $exp$ is homogeneous:
		$h_zp = exp(h_zq,\dot{h}_zv)(1)$, on the other hand, $h_zp = h_z\circ exp(q,v)(1)$.
		
		Let $f:\mathcal{U^\epsilon}\rightarrow M$, $p\mapsto exp(\pi(p), L_{\pi(p)}v)(1)$. Then $f_{|X} = id_X$, $df_p = L_p$ for every $p\in X$ and it is homogeneous.
	
	\noindent $\square$ 
	
	We have now all the tools we need to prove the following normal form results. 
	
	\begin{thm}[Homogeneous Weinstein Lagrangian neighborhood theorem]\label{thm: hWL}
		
		Let $(M,h)$ be a $2n$-dimensional $\mathbb{R}^\times$-manifold, let $X$ be an $n$-dimensional $\mathbb{R}^\times$-submanifold with $i:X\xhookrightarrow{} M$ the inclusion map, and let $\omega_0,\ \omega_1$ be two $1$-homogeneous symplectic forms on $M$ such that $i^*\omega_0 = i^*\omega_1 = 0$ ($X$ is Lagrangian for both). Then, there exists homogeneous neighborhoods $\mathcal{U}_0$ and $\mathcal{U}_1$  of $X$ in $M$ and a homogeneous diffeomorphism $\varphi:\mathcal{U}_0\rightarrow \mathcal{U}_1$ such that $\varphi^*\omega_1 = \omega_0$ and the following diagram commutes
		\begin{center}
			\begin{tikzpicture}
				\node at (0,0) (a) {$\mathcal{U}_0$};
				\node at (1,-1) (b) {$X$};
				\node at (2,0) (c) {$\mathcal{U}_1$};
				\draw [->] (b) to node[left] {$i$} (a);
				\draw [->] (a) to node[above] {$\varphi$} (c);
				\draw [->] (b) to node[right] {$i$}(c);
			\end{tikzpicture}
		\end{center}
	\end{thm}
	
	\noindent \textit{Proof:}
		Let us choose a Riemannian metric $g$ on $M$ that is invariant with respect to $h$; at each $p\in M$, $g_p(\cdot,\cdot)$ is a positive-define inner product.
		Fix $p\in X$ and let $V = T_pX$ and $V^\perp$ be the orthogonal complement of $V$ in $T_pM$, relative to $g_p(\cdot,\cdot)$.
		
		Since $i^*\omega_0=i^*\omega_1=0$, the space $V$ is Lagrangian subspace of both $(T_pM,\omega_{0|p})$ and $(T_pM,\omega_{1|p})$. By symplectic linear algebra, we obtain from $V^\perp$ a homogeneous linear isomorphism $L_p:T_pM\rightarrow T_pM$, such that $L_{p|V} = id_{V}$, $L^*_p\omega_{1|p} = \omega_{0|p}$ and depends smoothly on $p$.
		The existence of such linear isomorphisms $L_p$ is granted by Propositions 8.2 and 8.3 of \cite{da2001lectures} and their homogeneity can also directly be checked {(it's a consequence of homogeneity of $g$, implying homogeneity of choices of orthogonal complements as input for Proposition 8.3).}
		
		By theorem \ref{thm: whitney extension}, there are a homogeneous neighborhood $\mathcal{U}$ of $X$ and a homogeneous embedding $f:\mathcal{U}\rightarrow M$ with $f_{|X} = id_X$ and $df_p = L_p$, for every $p\in X$. Hence, $(f^*\omega_1)_p = (df_p)^*\omega_{1|p} = L_p^*\omega_{1|p} = \omega_{0|p}$. Applying the homogeneous Darboux-Weinstein Theorem \ref{thm: Darboux-Weinstein} to $\omega_0$ and $f^*\omega_1$, we find a homogeneous neighborhood $\mathcal{U}_0$ of $X$ and a homogeneous embedding $\phi:\mathcal{U}_0\rightarrow \mathcal{N}$ such that $\phi_{|X} = id_X$ and $\phi^*(f^*\omega_1) = \omega_0$ on $\mathcal{U}_0$.
		
		Set $\mathcal{U}_1= f(\mathcal{U})$, and $\varphi = f\circ \phi.$ 
	
	\noindent $\square$ 
	
	{
		\begin{thm}[Homogeneous Weinstein tubular neighborhood theorem]\label{App-thm-WLtub}
			
			Let $(M,\omega, h)$ be a $2n$-dimensional homogeneous symplectic manifold, let $X$ be an homogeneous Lagrangian submanifold, with $i:X\xhookrightarrow{} M$ the inclusion map, and $i_0:X\xhookrightarrow{} T^*X$ the Lagrangian embedding as the zero section. Let $\omega_0$ be the canonical symplectic form on $T^*X$.
			
			Then, there exist homogeneous neighborhoods $\mathcal{U}_0$ of $X$ in $T^*X$ and $\mathcal{U}$  of $X$ in $M$ and a homogeneous diffeomorphism $\varphi:\mathcal{U}_0\rightarrow \mathcal{U}$ such that $\varphi^*\omega = \omega_0$ and the following diagram commutes
			\begin{center}
				\begin{tikzpicture}
					\node at (0,0) (a) {$\mathcal{U}_0$};
					\node at (1,-1) (b) {$X$};
					\node at (2,0) (c) {$\mathcal{U}$};
					\draw [->] (b) to node[left] {$i_0$} (a);
					\draw [->] (a) to node[above] {$\varphi$} (c);
					\draw [->] (b) to node[right] {$i$}(c);
				\end{tikzpicture}
			\end{center}
		\end{thm}
	}
	
	{
		\noindent \textit{Proof:}
			The proof will follow exactly as in the standard case \cite[Theorem 9.3]{da2001lectures}, by using the homogeneous versions of the Tubular neighborhood theorem \ref{thm: tubular nbh}, and of the Weinstein Lagrangian neighborhood theorem \ref{thm: hWL}. In order for that strategy to work, all we need is to ensure that the canonical form $\omega_0$ on $T^*X$ is homogeneous, and that the isomorphism $NX\cong T^*X$ induced by $\omega$ is a homogeneous map.
			
			Let $(x,t)$ be homogeneous coordinates on $M$. Then $\omega_0 = dx\wedge \xi_x + dt\wedge \xi_t$. Computing $T^*h_z^*\omega_0$ we get
			\[(T^*h_z)^*\omega_0 = dx\wedge (z\xi_x) + d(zt)\wedge \xi_t = z\omega_0, \]
			so $\omega_0$ is homogeneous.
			
			Set $\tilde{\Omega}:N_xX\rightarrow T^*_xX,\ [v]\mapsto \omega(v,\cdot)_{|_{T_xX}}$. Take any $[v]\in N_xX$ and $u\in T_xX$. We want to prove that $\tilde{\Omega}(Th_z[v])(Th_z\cdot u) = (T^*h_z\cdot\tilde{\Omega}([v]))(Th_z\cdot u).$
			Since $X$ is $\mathbb{R}^\times$-invariant, $Th_z$ maps $T_xX$ to $T_{h_zx}X$ and descends to $NX$ (if $\gamma$ is a curve in $X$ with $\gamma(0)=x$ and $\dot{\gamma}(0) = u\in T_xX$, then $h_z\circ \gamma$ is a curve in $h_z(X)=X$ passing through $h_z(x)$ with tangent vector $Th_z\cdot u\in T_{h_z(x)}X$
			). By definition
			\[\tilde{\Omega}(Th_z[v])(Th_z\cdot u) = \omega_{h_z x}(Th_zv,Th_zu) = (h_z^*\omega)(v,u) = z\omega(v,u) = z\tilde{\Omega}([v])(u).\]
			On the other hand, we know that \cite{BruceGrabowskaGrabowski2017} $T^*h_z\xi = z\xi \circ (Th_{z^{-1}})_{|Th_z(x)X}$. So, $(T^*h_z\xi)(Th_z\cdot u) = z\xi(Th_{z^{-1}}\cdot Th_z\cdot u) = z\xi(u)$. Set $\xi = \tilde{\Omega}([v])$ and both sides agree and the diagram commutes 
			\begin{center}
				\begin{tikzpicture}
					\node at (0,0) (a) {$N_xX$};
					\node at (0,-1.5) (b) {$N_{h_z(x)}X$};
					\node at (2,0) (c) {$T^*_xX$};
					\node at (2,-1.5) (d) {$T^*_{h_z(x)}X$};
					\draw [->] (a) to node[left] {$Th_z$} (b);
					\draw[->] (b) to node[above] {$\tilde{\Omega}_0$} (d);
					\draw [->] (a) to node[above] {$\tilde{\Omega}_0$} (c);
					\draw [->] (c) to node[right] {$T^*h_z$} (d);
				\end{tikzpicture}
			\end{center}
			
			Therefore, the isomorphism $NX\cong T^*X$ is a homogeneous map. $\square$
		
		\noindent $\square$ 
	}
	
	\section{Homogeneous Hamiltonian dynamics}\label{anexo: H-J}
	
	In this section, we will derive the Hamilton-Jacobi equation by following the same steps as in \cite{Bravetti_2017}, but in the homogeneous setting. Let us start with a homogeneous symplectic manifold $(\Sigma,\omega,h)$, that is, $\omega$ is a symplectic form such that $h_z^*\omega = z\omega$. We want to see how Hamiltonian mechanics works in this scenario and what the conditions are to have a similar result to classical mechanics.
	
	Let $H\in C^\infty(\Sigma)$ be a function and define $X_H\in \mathfrak{X}(\Sigma)$ such that
	\[\omega(X_H,\cdot) = dH.\]
	
	Suppose $H$ is $1$-homogeneous, then $dH$ is also $1$-homogeneous form. We want to prove that the Hamiltonian vector field $X_H$ is $0$-homogeneous, that is, $(h_z)_*X_H = X_H$. So, using the relation 
	\begin{equation}\label{eq: relation}
		i_Xf^*\omega = f^*(i_{f_*X}\omega),
	\end{equation} let us compute $\omega((h_z)_*X_H,\cdot)$
	\begin{align*}
		z^{-1}i_{(h_z)_*X_H}\omega = i_{(h_z)_*X_H}(h_{z^{-1}})^*\omega =(h_{z^{-1}})^*i_{X_H}\omega = (h_{z^{-1}})^*dH = z^{-1}dH.
	\end{align*}
	So, $\omega((h_z)_*X_H,\cdot) = dH $, which by the definition of $X_H$ implies that $(h_z)_*X_H = X_H$. Since $X_H$ is $0$-homogeneous, its flow $\phi^s_{X_H}$ is also homogeneous, $\phi^s_{X_H}\circ h_z = h_z\circ \phi^s_{X_H}$.

	We want to see how Hamilton-Jacobi equation fits in homogeneity. Consider the homogeneous symplectic manifold $(\Sigma,\omega,h)$. Suppose that $(q^i,p_j),\ i,j=1,\dots,n$ are {homogeneous} Darboux coordinates on $\Sigma$, $\omega = \omega_{can}$ meaning that $h$ acts as
	\[h_z(q^1,\dots,q^n,p_1,\dots, p_n) = (q^1,\dots,q^n,zp_1,\dots, zp_n).\]
	This is possible if the homogeneous symplectic manifold is the cotangent bundle of another $\mathbb{R}^\times$-manifold, with canonical $1$-form $\alpha = p_idq^i$, that is also $1$-homogeneous.
	
	Consider the Hamiltonian $1$-homogeneous function $H\in C^\infty(\Sigma)$ and the respective $0$-homogeneous Hamiltonian vector field $X_H =\frac{\partial H}{\partial p_i}\frac{\partial}{\partial q^i} - \frac{\partial H}{\partial q^i}\frac{\partial}{\partial p_i}$. We say that a function is a \textit{homogeneous canonical transformation} if it preserves $\omega$ and is homogeneous.
	
	Suppose that we have such a function that has the change of coordinates $(q^i,p_i)\rightarrow(Q^i,P_i)$. Suppose now that the independent variables are $(q^i, Q^j)$. To make this canonical transformation homogeneous, $h_z$ now acts on the base, i.e., $h_z(Q^i,P_j) = (zQ^i,P_j)$. This makes sense because in the trivial case we have the canonical transformation $Q = p$ and $P = -q$.
	
	The invariance of $\omega$ implies that $\{Q^i,Q^j\} = \{P_i,P_j\} =0$ and $\{Q^i,P_j\} = \delta^i_j$ and also means that $\alpha$ is invariant up to an exact 1-homogeneous differential. So, there exists a function $F_1$ 1-homogeneous such that
	\[p_idq^i = P_idQ^i + dF_1.\] This function $F_1$ is called a \textit{homogeneous generating function}. We know that $dF_1 = \frac{\partial F_1}{\partial q^i}dq^i + \frac{\partial F_1}{\partial Q^i}dQ^i$, so
	\[\left(p_i-\frac{\partial F_1}{\partial q^i}\right)dq^i - \left(P_i + \frac{\partial F_1}{\partial Q^i}\right)dQ^i = 0.\]
	So, the homogeneous canonical transformation obeys the relation
	\[p_i = \frac{\partial F_1}{\partial q^i} \textnormal{ and } P_i = -\frac{\partial F_1}{\partial Q^i}.\]

	We want to study time-dependent Hamiltonian systems. Consider the extended phase space as $\Sigma^E = \Sigma\times\mathbb{R}$ and $h_z$ acting on $\Sigma^E$ as $h_z(q^i,p_i,t) = (q^i,zp_i,t)$. Consider also the Poincar\'{e}-Cartan 1-form $\eta_{PC}= p_idq^i - Hdt$. This is a $1$-homogeneous 1-form. Note that even though $\Sigma^E$ is a contact manifold, its Hamiltonian mechanics is given by the Lagrangian framework, which is a variational formulation \cite{abraham2008foundations,libermann1987symplectic}. In this case, the conditions for $X_H^E$ be a contact Hamiltonian vector field are
	\begin{align*}
		d\eta_{PC}(X_H^E) =0\qquad \text{ and }\qquad i_{X^E_H}dt = 1.
	\end{align*}
	Through a direct computation, we end with the following relation
	\[d\eta_{PC}(X_H^E) = 0 \Leftrightarrow X_H^E = X_H + \frac{\partial}{\partial t},\]where $X_H$ is our initial Hamiltonian vector field in $\Sigma$.
	
	\begin{remark}
		This extended Hamiltonian vector field satisfies also
		\[\eta_{PC}(X^E_H) = p_i\frac{\partial H}{\partial p_i} - H \qquad \text{ and }\qquad \mathcal{L}_{X_H^E} H = \frac{\partial H}{\partial t}.\]
	\end{remark}
	
	And since $X_H$ is 0-homogeneous, $X_H^E$ is also $0$-homogeneous. We need to find a homogeneous canonical transformation that leaves $d\eta_{PC}$ unchanged. So, if we add a $1$-homogeneous differential, it does not affect the equation. So,
	\[p_idq^i + Hdt -(P_idQ^i + Kdt) = dF_1,\] where $K$ is the new Hamiltonian. Choosing $(q^i,Q^j,t)$ as independent coordinates,
	\[\left(p_i-\frac{\partial F_1}{\partial q^i}\right)dq^i - \left(P_i+\frac{\partial F_1}{\partial Q^i}\right)dQ^i + \left(-K + H - \frac{\partial F_1}{\partial t}\right)dt = 0\]
	which implies that the 1-homogeneous generating function $F_1(q^i,Q^i,t)$ satisfies \[p_i = \frac{\partial F_1}{\partial q^i},\quad P_i = -\frac{\partial F_1}{\partial Q^i},\quad K = H - \frac{\partial F_1}{\partial t}.\]
	
	Now, we want the 1-homogeneous generating function $F_1$ such that the new Hamiltonian $K = 0$. Denote $F_1$ as $\mathbf{S}_t(q,Q,t)$. So, $\mathbf{S}_t$ satisfies the homogeneous Hamilton-Jacobi equation
	\begin{equation}\label{eq: H-J}
		H\left(q_i,\frac{\partial\mathbf{S}_t}{\partial q^i},t\right) = \frac{\partial\mathbf{S}_t}{\partial t}.
	\end{equation}
	

\end{document}